\documentclass{amsart}

\usepackage{microtype}
\usepackage{amssymb}
\usepackage{mathtools}
\usepackage{tikz-cd}
\usepackage{mathbbol} 

\usepackage[
	backend=biber,
	style=alphabetic, 
	backref=true,
	url=false,
	doi=false,
	isbn=false,
	eprint=false]{biblatex}

\renewbibmacro{in:}{} 
\AtEveryBibitem{\clearfield{pages}} 

\DeclareFieldFormat{title}{\myhref{\mkbibemph{#1}}}
\DeclareFieldFormat
[article,inbook,incollection,inproceedings,patent,thesis,unpublished]
{title}{\myhref{\mkbibquote{#1\isdot}}}

\newcommand{\doiorurl}{%
	\iffieldundef{url}
	{\iffieldundef{eprint}
		{}
		{http://arxiv.org/abs/\strfield{eprint}}}
	{\strfield{url}}%
}

\newcommand{\myhref}[1]{%
	\ifboolexpr{%
		test {\ifhyperref}
		and
		not test {\iftoggle{bbx:eprint}}
		and
		not test {\iftoggle{bbx:url}}
	}
	{\href{\doiorurl}{#1}}
	{#1}%
}

\usepackage[
	bookmarks=true,
	linktocpage=true,
	bookmarksnumbered=true,
	breaklinks=true,
	pdfstartview=FitH,
	hyperfigures=false,
	plainpages=false,
	naturalnames=true,
	colorlinks=true,
	pagebackref=false,
	pdfpagelabels]{hyperref}

\hypersetup{
	colorlinks,
	citecolor=blue,
	filecolor=blue,
	linkcolor=blue,
	urlcolor=blue
}

\usepackage[capitalize, noabbrev]{cleveref}
\crefname{subsection}{\S\!\!}{subsections}

\calclayout

\makeatletter
\@namedef{subjclassname@2020}{%
	\textup{2020} Mathematics Subject Classification}
\makeatother

\makeatletter
\renewcommand\subsubsection{\@startsection{subsubsection}{3}%
	\z@{.5\linespacing\@plus.7\linespacing}{0em}%
	{\normalfont\itshape}}
\makeatother
\newtheorem{theorem}{Theorem}
\newtheorem*{theorem*}{Theorem}
\newtheorem{proposition}[theorem]{Proposition}
\newtheorem*{proposition*}{Proposition}
\newtheorem{lemma}[theorem]{Lemma}
\newtheorem*{lemma*}{Lemma}
\newtheorem{corollary}[theorem]{Corollary}
\newtheorem*{corollary*}{Corollary}
\theoremstyle{definition}

\newtheorem*{definition*}{Definition}

\newtheorem*{remark*}{Remark}

\newtheorem*{example*}{Example}

\newtheorem*{construction*}{Construction}

\newtheorem*{convention*}{Convention}

\newtheorem*{terminology*}{Terminology}

\newtheorem*{notation*}{Notation}

\newtheorem*{question*}{Question}

\newcommand{\ot}{\otimes}

\newcommand{\N}{\mathbb{N}}
\newcommand{\Z}{\mathbb{Z}}

\newcommand{\R}{\mathbb{R}}
\renewcommand{\k}{\Bbbk}

\newcommand{\cyc}{\mathrm{C}}
\newcommand{\Ftwo}{{\mathbb{F}_2}}
\newcommand{\Fp}{{\mathbb{F}_p}}

\newcommand{\Fun}{\mathsf{Fun}}

\newcommand{\Top}{\mathsf{Top}}

\DeclarePairedDelimiter\bars{\lvert}{\rvert}

\DeclarePairedDelimiter\set{\{}{\}}

\newcommand{\id}{\mathsf{id}}
\renewcommand{\th}{\mathrm{th}}
\newcommand{\op}{\mathrm{op}}

\newcommand{\defeq}{\stackrel{\mathrm{def}}{=}}

\newcommand{\rC}{\mathrm{C}}

\newcommand{\rH}{\mathrm{H}}
\newcommand{\rI}{\mathrm{I}}

\newcommand{\rL}{\mathrm{L}}

\newcommand{\rP}{\mathrm{P}}
\newcommand{\rQ}{\mathrm{Q}}

\newcommand{\rU}{\mathrm{U}}

\newcommand{\cA}{\mathcal{A}}

\newcommand{\cC}{\mathcal{C}}
\newcommand{\cD}{\mathcal{D}}

\newcommand{\cM}{\mathcal{M}}

\newcommand{\cO}{\mathcal{O}}

\newcommand{\cX}{\mathcal{X}}
\newcommand{\cY}{\mathcal{Y}}
\newcommand{\cZ}{\mathcal{Z}}

\newcommand{\bC}{\mathbb{C}}

\newcommand{\bF}{\mathbb{F}}

\newcommand{\bS}{\mathbb{S}}

\newcommand{\VS}{\mathbb{VS}}
\DeclareMathOperator{\sus}{\Sigma}
\renewcommand{\Vec}{\mathsf{Vec}}

\newcommand{\kk}{\mathbb{k}}
\newcommand{\st}{\mathrm{st}}

\definecolor{darkdgmcolor}{rgb}{0.0, 0.0, 0.8}
\definecolor{darkgreen}{rgb}{0.0, 0.8, 0.0}
\definecolor{purple}{RGB}{153,50,204}
\definecolor{dgmcolor}{RGB}{255,20,147}
\definecolor{barccolor}{RGB}{20,147,255}
\definecolor{myred}{RGB}{227,26,28}

\definecolor{darkblue}{rgb}{0,0,0.7} 
\newcommand{\darkblue}{\color{darkblue}} 
\newcommand{\defn}[1]{{\darkblue \emph{#1}}} 

\DeclareMathOperator{\rad}{rad}
\newcommand{\VR}{\mathrm{VR}}
\DeclareMathOperator{\opH}{H}

\newcommand{\Sq}{\mathrm{Sq}}

\renewcommand{\theorem}{\medskip\noindent\textit{\darkblue Theorem}.\ }
\renewcommand{\lemma}{\medskip\noindent\textit{\darkblue Lemma}.\ }
\renewcommand{\proposition}{\medskip\noindent\textit{\darkblue Proposition}.\ }

\theoremstyle{theorem}
\newtheorem{introtheorem}{Theorem}

\DeclareMathOperator{\barc}{Bar}
\DeclareMathOperator{\img}{Im}
\renewcommand{\ker}{\operatorname{Ker}}

\newcommand{\field}{\k}

\newcommand{\rp}{\mathbb{RP}}

\newcommand{\diam}{\operatorname{diam}}

\newcommand{\dgh}{d_{\mathrm{GH}}}

\newcommand{\Coordinate}[2]%
{ \coordinate (#1) at (#2);
}

\newcommand{\thetabarc}[1]{\barc\img_{\theta}^{\VR} (#1)}

\newcommand{\dhi}{d_{\mathrm{HI}}}
\newcommand{\di}{d_{\mathrm{I}}}

\newcommand{\cost}{\mathrm{cost}}

\newcommand{\Hbarc}[3][\Ftwo]{\barc \opH_{#2}^{\VR} (#3)}
\newcommand{\sqbarcl}[3]{\barc\img_{\Sq^{#1}_{#2}}^{\VR} (#3)}

\newcommand{\fillrad}{\mathrm{Rad}_{\mathrm{fill}}}
\newcommand{\crit}{\mathrm{Rad}_{\mathrm{htpy}}}
\newcommand{\firstdeath}[2]{\mathrm{Rad}_{#1}(#2)}
\usepackage{float}
\usepackage{amsfonts}
\usepackage{graphicx}
\usepackage[all]{xy}
\usepackage{yfonts}
\usepackage{tikz}
\usepackage{pgfplots}
\pgfplotsset{compat=1.18}
\usepackage{etoolbox}
\usepackage{enumerate}
\usetikzlibrary{calc}
\usepackage{todo}
\title[Persistent cohomology operations]{Persistent cohomology operations \\ and Gromov--Hausdorff estimates}

\author{Anibal~M.~Medina-Mardones}
\address{Western University, Canada}
\email{\href{mailto:anibal.medina.mardones@uwo.ca}{anibal.medina.mardones@uwo.ca}}

\author{Ling~Zhou}
\address{Duke University, USA}
\email{\href{mailto:ling.zhou@duke.edu}{ling.zhou@duke.edu}}

\date{\today}


\begin{document}

\begin{abstract}
	This paper establishes the foundations of the theory of persistent cohomology operations, derives decomposition formulas for wedge sums and products, and proves their Gromov--Hausdorff stability.
	These results are used to construct pairs of Riemannian pseudomanifolds for which the Gromov--Hausdorff estimates derived from persistent cohomology operations are strictly sharper than those obtained using persistent homology.
\end{abstract}
	\maketitle
	\tableofcontents

\section{Introduction} \label{s:introduction}

When studying the shape of spaces equipped with a filtration, persistent homology has emerged as a crucial tool in both applied and theoretical topology.
By tracking the evolution of homological features across the filtration, this technique has found widespread applications in fields ranging from data analysis and machine learning (\cite{carlsson2013viral, lee2017quantifying}) to symplectic geometry and functional analysis (\cite{polterovich2020persistence, medina2023fuct_top}).
Despite its successes, persistent homology has limitations that mirror those of classical homology in the study of unfiltered topological spaces.
Just as homology fails to distinguish spaces with different homotopy types, persistent homology can overlook significant structural aspects of their filtered counterparts.

To address these shortcomings, researchers have begun exploring extensions of persistent homology that incorporate more refined homotopical information.
For example, the cup product in cohomology has been studied in \cite{contessoto_et_al:LIPIcs.SoCG.2022.31, memoli2024persistent, huang2005cup, yarmola2010persistence, contreras2022persistent} and some of its derived structures in \cite{herscovich2018higher, ginot2019distances, belchi2022persistence, zhou2023persistent, hess2024minimalmodels}.
In this work, we concentrate on cohomology operations, which are natural transformations between cohomology functors.

The first major developments integrating cohomology operations into the persistence framework were presented in \cite{medina2022per_st}; for completeness, we also mention earlier approaches with related goals in \cite{aubrey2011thesis} and \cite{postol2023persistence}.
Using new methods for the computation of Steenrod squares (\cite{medina2023fast_sq}), this work proposed and implemented algorithms for computing the \textit{Steenrod barcodes} of finite filtered simplicial complexes, with the resulting tool, \href{https://steenroder.github.io/steenroder/}{\texttt{steenroder}}, used to detect the nontrivial presence of these invariants on molecular data.

Despite the conceptual and computational advances of \cite{medina2022per_st}, a comprehensive exploration of persistent cohomology operations has remained absent, leaving key questions about their theoretical properties and broader applicability unanswered.

In this paper, we systematically develop a general and firm foundation for \textit{persistent cohomology operations} and present applications of these invariants to Riemannian geometry.

\medskip More precisely, let \(\k\) be a field.
A \textit{\(\k\)-linear cohomology operation} \(\theta\) is a natural transformation
\[
\theta \colon \rH^\ell(-;\k) \to \rH^m(-;\k),
\]
between cohomology functors.
From now on we will omit \(\k\) from the notation.

Given a functor \(X\) from the poset category $\R$ of real numbers to the category of cellular spaces, the image and kernel of \(\theta\) define two persistent modules, denoted \(\img_\theta(X)\) and \(\ker_\theta(X)\) respectively.
In this introduction we focus solely on the first of these, mentioning that all statements also hold for the second.

When \(X\) is the Vietoris--Rips filtration \(\VR(\cX)\) of a metric space \(\cX\), we use the simplified notation \(\rH_m^\VR(\cX)\) and \(\img_\theta^\VR(\cX)\) instead of \(\rH_m(\VR(\cX))\) and \(\img_\theta(\VR(\cX))\).

\medskip\noindent\textsc{Decomposition formulas.}
Our first contributions are the following theorems, proven respectively in \cref{ss:wedge sum} and \cref{ss:products}, providing decomposition formulas for the persistent cohomology modules of wedge sums and products of metric spaces.

\begin{introtheorem}\label{thm:decomposition1}
	For any two pointed metric spaces $\cX$ and $\cY$, a linear cohomology operation \(\theta\) satisfies:
	\[
	\img_\theta^\VR(\cX \vee \cY) \cong \img_\theta^\VR(\cX) \oplus \img_\theta^\VR(\cY),
	\]
	where the wedge sum \(\cX \vee \cY\) is equipped with the gluing metric.
\end{introtheorem}

\begin{introtheorem}\label{thm:decomposition2}
	For any two metric spaces $\cX$ and $\cY$, over any prime, the total Steenrod operation \(\rP\) satisfies:
	\[
	\img^{\VR}_\rP(\cX \times \cY) \cong \img^{\VR}_\rP(\cX) \ot \img^{\VR}_\rP(\cY),
	\]
	where the product \(\cX \times \cY\) is equipped with the \(\ell^\infty\) metric.
\end{introtheorem}

These follow from corresponding statements for wedge sums and products of general cellular $\mathbb{R}$-spaces presented, respectively, in \cref{ss:wedge_sum_general} and \cref{ss:product_general}.

\medskip\noindent\textsc{Gromov--Hausdorff Inequalities.}
A fundamental property of persistent homology is its stability, which has multiple formulations; here, we focus on the Vietoris--Rips stability for metric spaces, as presented in the following.
\begin{theorem*}[\cite{chazal2009gromov, chazal2014geometric}]
	\label{thm:stability ph}
	For any two metric spaces $\cX$ and $\cY$ and non-negative integer $m$:
	\begin{equation}\label{eq:gh_inequality_intro_homology}
		\di\big(\rH_m^\VR(\cX),\, \rH_m^\VR(\cY)\big) \leq 2 \cdot \dgh(\cX,\cY),
	\end{equation}
	where \(\di\) and \(\dgh\) are respectively the interleaving and extended Gromov--Hausdorff distances.
\end{theorem*}

Our next contribution is to show that persistent cohomology operations are also stable in this sense.

More precisely, in \cref{thm:stability VR}, we prove the following result using the homotopy interleaving distance of \cite{blumberg2023interleaving}.

\begin{introtheorem}\label{thm:stability intro}
	For any two metric spaces $\cX$ and $\cY$ and a linear cohomology operation~$\theta$:
	\begin{equation}\label{eq:gh_inequality_intro}
		\di\big(\img_\theta^\VR(\cX),\, \img_\theta^\VR(\cY)\big) \leq 2 \cdot \dgh(\cX,\cY).
	\end{equation}
\end{introtheorem}

This follows from a statement holding for general cellular $\mathbb{R}$-spaces presented in \cref{thm:stability theta}.

\medskip Given the inequalities \eqref{eq:gh_inequality_intro_homology} and \eqref{eq:gh_inequality_intro},
it is natural to ask whether persistent cohomology operations can yield sharper lower bounds for the Gromov--Hausdorff distance than those obtained via persistent homology.
This would demonstrate the greater discriminating power of persistent cohomology operations compared to persistent homology.

We employ certain \textit{critical radii}, inspired by Gromov's filling radius~\cite{gromov1983filling}, to establish this result in \cref{ss:distance_estimate_rpn}.
Specifically, we consider the real projective space \(\rp^n\), with diameter \(\pi\), obtained as the quotient of a round sphere under the antipodal action, and the wedge sum of round spheres \(\bS_{\mathbb{RP}^n} = \bS^1 \vee \dots \vee \bS^n\), each also with diameter \(\pi\).

\begin{introtheorem}\label{thm:inequalities intro}
	For every \(n > 1\) the following inequalities hold:
	\begin{enumerate}
		\item For any \(m \in \N\)
		\[
		\di\Big(\rH^\VR_m(\rp^n),\, \rH^\VR_m(\bS_{\rp^n})\Big) < \frac{\pi}{4}.
		\]

		\item There is \(k \in \N\) such that
		\[
		\di\Big(\img_{\Sq^k}^\VR(\rp^n),\, \img_{\Sq^k}^\VR(\bS_{\rp^n})\Big) \geq \frac{\pi}{3},
		\]
		where \(\Sq^k\) is the \(k^\th\) Steenrod square.
	\end{enumerate}
\end{introtheorem}

\medskip As explained in \cref{ss:genberal_distance_comparison}, our arguments extend to Lens spaces and odd prime Steenrod operations, provided certain concrete statements about their critical radii are established.
Investigating these geometric questions falls outside the scope of this work.

\section*{Related work}

\subsection*{Decomposition formulas}

Both \cref{thm:decomposition1} and \cref{thm:decomposition2} build on previous studies of the Vietoris--Rips filtration of metric spaces.
Specifically of their wedge sums \cite{adamaszek2020homotopy} and products \cite{adamaszek2017VietorisProduct, gakhar2019k, lim2024vietoris}.

\subsection*{Stability}

Our proof of the stability of persistent cohomology operations builds upon the strategies developed in \cite{blumberg2023interleaving, zhou2023beyond, memoli2024persistent}.

In \cite[Theorem~87]{ginot2019distances}, Ginot and Leray prove a Gromov--Hausdorff inequality using the \textit{Steenrod interleaving distance}.
This is the interleaving distance in the category of persistent algebras that are also persistent modules over the Steenrod algebra \(\mathcal{A}_p\).
Our stability result, \cref{thm:stability intro}, holds for all cohomology operations.
In the important case of Steenrod operations, it seems likely that our result could be deduced from theirs.

\subsection*{Gromov--Hausdorff estimates via interleaving estimates}

The exact value of the Gromov--Hausdorff distance is known in only a small number of cases \cite{memoli2012some,ji2021gromov,adams2022gromov,talipov2022gromov,lim2021gromov,harrison2023quantitative,saul2024gromov,saul2024some}.
Much of the effort in these works, as well as in the broader study of the Gromov--Hausdorff distance, has been devoted to establishing effective bounds on this distance (see \cite{lim2021gromov} for a comprehensive discussion of the literature).

In this paper, we contribute to the line of work about estimating Gromov--Hausdorff distances between metric spaces by seeking lower bounds based on persistent cohomology operations.
Indeed, \cref{thm:inequalities intro} illustrates the stronger discriminating power of persistent cohomology operations over (standard) persistent homology.
It does so by providing bounds on the interleaving distance of the persistent homology and persistent cohomology operations of the Vietoris--Rips filtration of pairs of Riemannian pseudomanifolds.
For example, combining \cref{thm:stability intro} and \cref{thm:inequalities intro}, we establish a lower bound of $\frac{\pi}{6}$ for the Gromov--Hausdorff distances between \(\rp^n\) and \(\bS_{\rp^n}\), which is strictly greater than the lower bound obtained from persistent homology, which does not exceed $\frac{\pi}{8}$.
See \cref{s:gh_estimates} for more details.

This approach, which employs topological invariants more informative than persistent homology, builds upon prior work explored in \cite{zhou2023beyond, zhou2023persistent, memoli2024persistenthomotopy, memoli2024persistent, memoli2025ephemeral}.
These studies have investigated sharper lower bounds on the Gromov--Hausdorff distance by leveraging persistent invariants derived from homotopy groups, rational homotopy groups, the cohomology ring, the Lyusternik--Schnirelmann category, Sullivan minimal models, and filtered chain complexes.

\subsection*{Critical radii}

The Kuratowski embedding \(K \colon \cX \to \rL^\infty(\cX)\) of a metric space \(\cX\) is its canonical isometric inclusion, given by \(\cX \ni x \mapsto d_\cX(x,\cdot)\), into the Banach space $\rL^\infty(\cX)$.
For any $r>0$, let \(\rU_r(\cX)\) denote the radius \(r\) neighborhood of \(K(\cX)\) in \(\rL^\infty(\cX)\).
It was established in \cite{chazal2009gromov} for finite metric spaces and later extended to compact metric spaces in \cite{lim2024vietoris} that the Vietoris--Rips complex \(\VR_{2r}(\cX)\) is naturally homotopy equivalent to \(\rU_r(\cX)\) for every \(r > 0\).

Gromov's filling radius is defined as the smallest \(r\) for which the fundamental class of a closed connected Riemannian manifold becomes null-homologous in \(\rU_r(\cM)\).
The use of this invariant and other similar critical radii is crucial in our proof of \cref{thm:inequalities intro}.
This approach is similar to the one used in \cite[Proposition~9.38]{lim2024vietoris} where the authors consider pairs of spaces consisting of an $n$-manifold \(\cM\) and a space \(\cX\) with trivial \(n\)-homology, and then bound the interleaving distance of the \(n\)-th persistent homology by the filling radius of \(\cM\).
In contrast, we consider a space \(\cX\) whose homology matches that of \(\cM\) in all degrees.

\subsection*{Recent developments}

After the posting of the present work, F.~M\'emoli and the second-named author studied Gromov--Hausdorff estimates for \(\rp^n\) and \(\bS^n\) using persistent invariants derived from cup-length, the Lyusternik--Schnirelmann category, topological complexity, and zero-divisor cup-length \cite{memoli2025persistence}.
Their analysis also yields the lower bound \(\frac{\pi}{3}\) obtained here via persistent cohomology operations.


\section*{Acknowledgements}

We thank Henry Adams, Liam Barham, Gr\'egory Ginot, Kathryn Hess, Umberto Lupo, and Dennis Sullivan for insightful comments.
We are especially grateful to Facundo M\'{e}moli for his thoughtful insights and valuable input on key aspects of the relevant literature.
We also thank the reviewers for their careful reading and many suggestions, which significantly improved the clarity and presentation of this work.

A.M. gratefully acknowledges the excellent working conditions provided by the Max Planck Institute for Mathematics in Bonn and the research unit LAGA at Université Sorbonne Paris Nord, where part of this work was carried out.

L.Z. gratefully acknowledges the support of The Ohio State University and the Institute for Computational and Experimental Research in Mathematics at Brown University, where part of this work was carried out.

A.M. received support for this project from grants ANR 20 CE40 0016 01 PROJET HighAGT and NSERC Discovery 06445-2024 RGPIN.
L.Z. received support from an AMS--Simons Travel Grant.


\section{Persistence and Vietoris--Rips filtrations}

In this section we provide an overview of the key concepts in the theory of persistence and its applications to metric geometry.

\subsection{$\R$-diagrams and interleaving distance}

\subsubsection{}

We think of the set of real numbers $\R$ as a category via its poset structure.

For any category $\cC$, the category $\cC^\R$ of \defn{$\R$-diagrams in $\cC$} consists of functors from $\R$ to $\cC$ and natural transformations between them.

For any $\R$-diagram $F$, we denote the image of $r \in \R$ by $F_r$ and the image of a morphism $s \leq t$ by $F_{s,t}$.

\subsubsection{}\label{ss:interleaving}

Let $\delta$ be a fixed non-negative real number.
For any $\R$-diagram $F \colon \R \to \cC$, we write $F^\delta$ for the $\R$-diagram defined by
\[
F^\delta_r = F_{r+\delta}
\quad\text{and}\quad
F^\delta_{s \leq t} = F_{s+\delta \leq t+\delta}.
\]

A \defn{$\delta$-interleaving} of two $\R$-diagrams $F,G \colon \R \to \cC$ is a pair of natural transformations
$\phi \colon F \to G^\delta$ and $\psi \colon G \to F^\delta$ such that for any $r \in \R$ the morphism $\psi_{r+\delta} \circ \phi_r$ and $\phi_{r+\delta} \circ \psi_r$ are respectively equal to $F_{r,r+2\delta}$ and $G_{r,r+2\delta}$.

The \defn{interleaving distance} between them is defined by
\[
\di(F, G) = \inf\set{\delta \geq 0 \mid F \text{ and } G \text{ are } \delta\text{-interleaved}}.
\]

Functorially changing the target category $\cC$ does not increase the interleaving distance.

\proposition (\cite[page~1508]{bubenik2015metrics})
For any functor $\Phi \colon \cC \to \cD$ we have
\[
\di(\Phi \circ F, \Phi \circ G) \leq \di(F, G)
\]
for any pair of functors $F, G \colon \R \to \cC$.

\medskip\noindent
This is because the $\delta$-interleaving functors of $F$ and $G$ compose with $\Phi$ to produce a $\delta$-interleaving of $\Phi \circ F$ and $\Phi \circ G$.

\subsection{$\R$-spaces and homotopy interleaving distance}

\subsubsection{}\label{ss:R-spaces}

Let $\Top$ denote the category of compactly generated weak Hausdorff spaces and continuous maps.
We write $\Top^\R \defeq \Fun(\R,\Top)$ for the category of \defn{$\R$-spaces}--where $\R$ is regarded as a poset category--equipped with the projective model structure in which weak equivalences and fibrations are defined objectwise.
An $\R$-space $X$ is \defn{cellular} if each $X_r$ is a CW-complex.

We say that a morphism $f \colon X \to Y$ in $\Top^\R$ is a \defn{weak equivalence} if, for every $r \in \R$, the map $f_r \colon X_r \to Y_r$ is a weak homotopy equivalence.
Two $\R$-spaces $X$ and $Y$ are \defn{weakly equivalent}, denoted $X \simeq Y$, if they become isomorphic in the homotopy category or, equivalently, $X \simeq Y$ if there exists an $\R$-space $W$ and weak equivalences
\[
X \xleftarrow{\ \sim\ } W \xrightarrow{\ \sim\ } Y
\]
\cite[\S2.4]{blumberg2023interleaving}.

\subsubsection{}\label{def:dhi}

Two $\R$-spaces $X$ and $Y$ are said to be \defn{$\delta$-homotopy-interleaved} for some $\delta \geq 0$, if there exist $\R$-spaces $X' \simeq X$ and $Y' \simeq Y$ such that $X'$ and $Y'$ are $\delta$-interleaved.

Following \cite{blumberg2023interleaving}, the \defn{homotopy interleaving distance} between two $\R$-spaces $X$ and $Y$ is given by
\[
\dhi(X,Y) = \inf \set{\delta \geq 0 \mid X,Y \text{ are }\delta\text{-homotopy-interleaved}}.
\]

\subsection{Persistence modules and bottleneck distance}

\subsubsection{}

Let $\Vec$ be the category of vector spaces and linear maps over a fixed field $\k$.
The category $\Vec^\R$ is referred to as that of \defn{persistence modules}.

We denote, with slight abuse of notation, by $\kk$ and $0$ the functors that assign $\kk$ and $0$, respectively, to every object of $\R$ and the identity map to each morphism.

\subsubsection{}

Let $I$ be an interval in $\R$.
The persistence module $\kk[I]$, referred to as an \defn{interval module}, is defined explicitly by
\[
\kk[I]_t =
\begin{cases}
	\kk & \text{if } t \in I, \\
	\hfil 0 & \text{otherwise},
\end{cases}
\qquad \qquad
\kk[I]_{s, t} =
\begin{cases}
	\id & \text{if } s, t \in I, \\
	\hfil 0 & \text{otherwise}.
\end{cases}
\]

A persistence module $V$ is said to be \defn{interval decomposable} if there is a multiset of intervals $(I_\lambda)_{\lambda \in \Lambda}$ such that $V$ is isomorphic to $\bigoplus_{\lambda \in \Lambda} \kk[I_\lambda]$.
In this case we refer to the multiset as the \defn{barcode} of $V$ and denote it $\barc V$.
By Azumaya’s theorem \cite{azumaya1950theorem}, $\barc V$ is unique up to reordering.

\subsubsection{}

Crawley-Boevey's result \cite{Crawley-Boevey.2015} is a widely-used existence theorem for barcode decompositions.
It ensures that any persistence module \(V\) with \(V_t\) being a finite dimensional vector space for every \(t \in \R\) has a barcode decomposition.
However, this condition is often too restrictive.

A more general condition is q-tameness.
For a detailed presentation of the approach we summarize below, consult \cite{Chazal.2016a, Chazal.2016b}.
A persistence module \(V\) is \defn{q-tame} if the rank of the map \(V_{s,t} \colon V_s \to V_t\) is finite for all \(s < t\).

For an interval \(I \subset \R\), let \(C(I)\) denote the interval module, ie the persistence module that is \(\Bbbk\) on \(I\), zero elsewhere, and whose structure maps are the identity within \(I\) and zero otherwise.
The infinite product \(\prod_{n \in \N_{> 0}} C([0,1/n))\) demonstrates that not all q-tame persistence modules admit a barcode decomposition in the traditional sense.
A persistence module \(V\) is termed \defn{ephemeral} if all the maps \(V_{s,t} \colon V_s \to V_t\) are zero for \(s < t\).

The \defn{radical} \(\rad V\) of a persistence module \(V\) is defined as the unique minimal submodule of \(V\) such that the cokernel of the inclusion \(\rad V \hookrightarrow V\) is an ephemeral persistence module.
Specifically, \((\rad V)_t = \sum_{s<t} \img V_{s,t}\).
For example, the radical of the infinite product \(\prod_{n \in \N_{> 0}} C([0, 1/n))\) is the direct sum \(\bigoplus_{n \in \N_{> 0}} C((0,1/n))\).
If \(V\) is q-tame its radical has a barcode decomposition, which describes the isomorphism type of \(V\) ``up to ephemerals".\footnote{
This idea is formalized using the so-called \textit{observable category}, which is equivalent to the quotient of the category of persistence modules by the subcategory of ephemeral persistence modules.
The barcode of the radical of a q-tame persistence module \(V\) serves as a complete invariant of \(V\) in the observable category.}
We define the barcode of a q-tame persistent module as the barcode of its radical.

\subsubsection{}
\label{ss:algebraic_stability}
Let $A$ and $A'$ be two possibly empty multisets of points in $\R^2$ lying above the diagonal.
A subset $P \subset A \times A'$ is said to be a \defn{partial matching} between $A$ and $A'$ if every point $(a, b) \in A$ is matched with at most one point of $A'$ and every point $(a', b') \in A'$ is matched with at most one point of $A$.

The \defn{bottleneck distance} between $A$ and $A'$ is defined as
\[
d_B(A, A') = \
\inf \set[\big]{\mathrm{cost}(P) \mid P \subset A \times A' \text{ is a partial matching}}
\]
where $\mathrm{cost}(P)$ is the largest between
\[
\sup \set[\Big]{\max\set{|a - a'|, |b - b'|} \bigm| (a, b) \in A, (a', b') \in A' \text{ are matched in }P}
\]
and
\[
\sup \set[\Big]{\tfrac{ |a - b|}{2} \bigm| (a, b) \in A \cup A' \text{ is unmatched}}.
\]

\proposition (\cite[Theorem~5.14]{Chazal.2016a})
The interleaving distance of q-tame persistence modules agrees with the bottleneck distance of its barcodes.

\subsubsection{}

For any integer $k \geq 0$, applying the degree $k$ (reduced) homology functor (with coefficients in~$\kk$) to an $\R$-space produces its \defn{persistent $k$-homology}, a persistence module denoted $\rH_k(X; \kk)$.
We will also consider the reduced version of this construction.

It is also convenient to consider functors from $\R^\op$ to $\Vec$, which, by abuse of terminology, we also refer to as persistence modules.
The most important example of these is the \defn{persistent $k$-cohomology} $\rH^k(X; \kk)$ of an $\R$-space $X$.

For simplicity, we often omit the coefficient field $\kk$ when the results do not depend on the choice of coefficient fields.

An $\R$-space is said to be \defn{q-tame} if its persistent homology and cohomology are both q-tame as persistence modules for any degree and field of coefficients.

\subsection{Vietoris--Rips filtrations and stability}

\subsubsection{}

Given a metric space $\cX$ and $r > 0$, let $\VR_r(\cX)$ denote the geometric realization of the simplicial complex whose vertices are the points of $\cX$ and whose simplices are the finite subsets of $\cX$ with diameter strictly less than $r$.
If $s \leq t$ then $\VR_s(\cX)$ is a subset of $\VR_t(\cX)$.
So, together with the inclusion maps, these spaces define an $\R$-space with $\VR_r(\cX) = \emptyset$ for non-positive values of $r$.
We refer to it as the \defn{Vietoris--Rips filtration} of $\cX$.
For any \(m \in \N\) we denote the persistence module \(\rH_m\big(\VR(\cX)\big)\) by \(\rH_m^{\VR}(\cX)\).

\subsubsection{}\label{sss:q-tameness}

Recall that a metric space is said to be \defn{totally bounded} if, for every $\epsilon > 0$, there exists a finite set of points such that the entire space can be covered by the union of $\epsilon$-balls centered at these points.

\proposition (\cite[Proposition~5.1]{chazal2014geometric})
If $\cX$ is a totally bounded metric space, then $\VR(\cX)$ is $q$-tame.

\subsubsection{}

For notational simplicity, when referring to intervals (or bars) in a barcode of $\VR(\cX)$, we will use the notation $(a, b)$ to denote the left-open, right-closed interval $(a, b]$.
This should not cause confusion due to the following.

\proposition (\cite[Theorem~5.2]{lim2024vietoris})
If $\cX$ is compact, the barcode of the persistent homology of $\VR(\cX)$ consists of only intervals that are left-open and right-closed.

\subsubsection{}\label{thm:stability-HI}

The \defn{Gromov--Hausdorff distance} $\dgh(\cX, \cY)$ between two compact metric spaces $\cX$ and $\cY$ is defined as
\[
\dgh(\cX, \cY) = \inf_{\cZ, \varphi, \psi} \set[\Big]{d_\rH\big(\varphi(\cX), \psi(\cY)\big)},
\]
where the infimum is taken over all metric spaces $\cZ$ and isometric embeddings $\varphi \colon \cX \to \cZ$ and $\psi \colon \cY \to \cZ$, and $d_\rH$ denotes the Hausdorff distance.

\proposition (\cite[Theorems~1.9 and~1.10]{blumberg2023interleaving})
Let $\cX$ and $\cY$ be two metric spaces.
Then, for any $k \in \N,$
\[
\di\big(\rH_k^\VR(\cX), \rH_k^\VR(\cY)\big) \leq
\dhi\big(\VR(\cX), \VR(\cY)\big) \leq
2 \cdot \dgh(\cX, \cY).
\]


\section{Persistent cohomology operations}\label{s:steenrod}

The first part of this section provides a review of classical material on cohomology operations, following standard references such as \cite{steenrod1962cohomology, mosheroperations1968, may1970general}.

Subsequently, we introduce persistent invariants derived from these operations, establishing sum and product decompositions and proving a stability result for them.

\subsection{Eilenberg--MacLane spaces}

\subsubsection{}

Let $\pi$ be an abelian group.
We denote by $K(\pi, n)$ any connected cellular space which has only one non-trivial homotopy group, namely, $\pi_n(K(\pi, n)) = \pi$.
Such a space is unique up to homotopy, and is referred to as an \defn{Eilenberg--MacLane space}.
For example, since the circle has the real line as its universal cover and fundamental group isomorphic to $\Z$, it is a model for $K(\Z,1)$.

Let $\cX$ be a space with the homotopy type of a connected cellular space.
Unless otherwise indicated, we use reduced homology and cohomology.
The \defn{Hurewicz homomorphism} $h \colon \pi_i(\cX) \to \rH_i(\cX)$ is defined by choosing a generator $u$ of $\rH_i(\bS^i)$ and sending the homotopy class $[f]$ of a based map $f \colon \bS^i \to \cX$ to $f_*(u)$, where $\bS^i$ denotes the $i$-sphere.
Using the properties of this map and the universal coefficient theorem for cohomology one shows the existence of a preferred class $\iota_n \in \rH^n(K(\pi, n); \pi)$ termed the \defn{fundamental class} of $K(\pi, n)$.

The relevance of Eilenberg--MacLane spaces for the study of cohomology is the following bijection holding for any space \(\cX\) with the homotopy type of a CW complex:
\[
\begin{tikzcd}[column sep=small, row sep=0]
	{[\cX, K(\pi, n)]} \rar & \rH^n(\cX; \pi) \\
	{[f]} \rar[maps to] & f^*(\iota_n),
\end{tikzcd}
\]
where $[\cX, K(\pi, n)]$ denotes the set of homotopy classes of based maps from $\cX$ to $K(\pi, n)$.
We refer to this bijection as the \defn{representability of cohomology}.

\subsection{Cohomology operations}

\subsubsection{}
Let both $\pi$ and $G$ be abelian groups.
A \defn{cohomology operation} of type $(\pi, n; G, m)$ is a family of functions
\[
\theta_\cX \colon \rH^n(\cX; \pi) \to \rH^m(\cX; G),
\]
one for each cellular space $\cX$, satisfying the naturality condition $f^* \theta_{\cY} = \theta_{\cX} f^*$ for any map $f \colon \cX \to \cY$.
We will denote by $\cO(\pi, n; G,m)$ the set of cohomology operations of type $(\pi, n; G,m)$.

\subsubsection{}
There is a bijection between $\cO(\pi, n; G,m)$ and homotopy classes of based maps between their Eilenberg--MacLane spaces.
Explicitly, given $[f] \in [K(\pi,n), K(G,m)]$ and an $n$-cohomology class in a space $\cX$ with $\pi$ coefficients, say $g^*(\iota_n)$ for some $g \colon \cX \to K(\pi, n)$, the value of the cohomology operation $\theta_{\cX}^{[f]}$ on it is
\[
\theta_{\cX}^{[f]}(g^*(\iota_n)) = (f \circ g)^*(\iota_m).
\]
In other words, cohomology operations are parameterized by the cohomology of Eilenberg--MacLane spaces, ie,
\[
\cO(\pi, n; G,m) \cong \rH^m(K(\pi,n); G).
\]

\subsection{Steenrod operations}\label{ss:steenrod}

\subsubsection{}

The \defn{suspension} of a space $\cX$, denoted as $\sus \cX$, is the product of $\cX$ with the interval $[0,1]$ having each $\cX \times \set{0}$ and $\cX \times \set{1}$ collapsed to a point.
For example, $\sus \bS^n$ is homotopy equivalent to $\bS^{n+1}$.

Reduced cohomology satisfies a natural isomorphism
\[
\rH^n(\cX) \cong \rH^{n+1}(\sus \cX)
\]
with any coefficients, referred to as the \defn{suspension isomorphism}.

\subsubsection{}

A \defn{stable cohomology operation of degree \(k\)} is a natural family of cohomology operations
\[
\set{\theta^n \in \cO(\pi, n; G, n+k)}_{n \in \N}
\]
commuting with suspension isomorphisms, ie, making the diagram
\[
\begin{tikzcd}
	\rH^n(\cX; \pi) \rar["\theta^n_\cX"] \dar["\cong" left] & \rH^{n+k}(\cX; G) \dar["\cong" right] \\
	\rH^{n+1}(\sus\cX; \pi) \rar["\theta^{n+1}_{\sus\cX}"] & \rH^{n+1+k}(\sus\cX; G)
\end{tikzcd}
\]
commute for each \(\cX\) and \(n\).
We denote by \(\cO^\st_k(\pi, G)\) the set of stable cohomology operations of degree \(k\).

Using the representability of cohomology, \(\cO^\st_k(\pi, G)\) can be described as the inverse limit
\[
\cO^\st_k(\pi, G) \cong \varprojlim_n \rH^{n+k}(K(\pi, n); G).
\]
Explicitly, a chosen map
\[
s_n \colon \sus K(\pi, n) \to K(\pi, n+1)
\]
adjoint to a weak equivalence \(K(\pi, n) \simeq \Omega K(\pi, n+1)\) induces the transition map
\[
\rH^{n+1+k}(K(\pi, n+1); G) \xrightarrow{\ s_n^*\ } \rH^{n+1+k}(\sus K(\pi, n); G) \xrightarrow{\ \cong\ } \rH^{n+k}(K(\pi, n); G),
\]
where the second map is the inverse of the suspension isomorphism.

\subsubsection{}

We will be interested in combining all stable operations over a fixed coefficient field $\k$ to form, via compositions, the (graded) algebra:
\[
\cO^\st(\k) = \bigoplus_k \cO^\st_k(\k, \k).
\]

For any prime $p$ let \(\Fp\) be the field with \(p\) elements.
The algebra $\cO^\st(\Fp)$ is denoted $\cA_p$ and referred to as the mod~$p$ \defn{Steenrod Algebra}.

\subsubsection{}

The Steenrod algebra $\cA_2$ is generated by the \defn{Steenrod squares} $\set{\Sq^i}_{i \in \N}$, where $\Sq^i \colon \rH^n(\cX; \Ftwo) \to \rH^{n+i}(\cX; \Ftwo)$ for any $n\geq 0$.
These are degree $i$ stable cohomology operations which are axiomatically characterized by the following properties, holding for all cohomology classes $\alpha,\beta \in \rH^*(\cX; \Ftwo)$ in any space $\cX$:
\begin{enumerate}
	\item \(\Sq^0(\alpha) = \alpha,\)
	\item \(\Sq^i(\alpha) = 0, \quad \text{if } i > \deg\alpha,\)
	\item \(\Sq^i(\alpha) = \alpha^2, \quad \text{if } i = \deg\alpha,\)
	\item \(\Sq^i(\alpha \beta) = \textstyle\sum_{j=0}^{i} \Sq^j(\alpha) \Sq^{i-j}(\beta).\)
\end{enumerate}

The \defn{total Steenrod square} is the (inhomogeneous) cohomology operation
\[
\Sq = \Sq^0 + \Sq^1 + \Sq^2 + \dotsb
\]
acting on \(\bigoplus_{m \in \N} \rH^m(\cX; \Ftwo)\).

We mention that $\Sq^1$ is the Bockstein homomorphism of the exact sequence
\[
0 \to \Z/2 \to \Z/4 \to \Z/2 \to 0.
\]

\subsubsection{}\label{sss:cohomology_rpn}

In this and the following examples we use unreduced cohomology.

Consider the system of unit spheres
\[
\bS^0 \subset \bS^1 \subset \bS^2 \subset \dotsb,
\]
defined by the canonical inclusions \(\R^n \subset \R^{n+1}\), with the \(\cyc_2\) action defined by the diagonal action of \(\set{-1,1} \subset \R\).
A model for \(K(\Z/2,1)\) is the infinite real projective space \(\rp^\infty\), defined as the union of the sequence of (finite) \defn{real projective spaces}
\[
\rp^0 \subset \rp^1 \subset \rp^2 \subset \dotsb,
\]
where each \(\rp^n\) is the space of orbits of \(\bS^n\) under the \(\cyc_2\)-action.

The cohomology algebra of $\rp^\infty$ with mod 2 coefficients is the polynomial algebra generated by a single element $\sigma$ in degree 1.
Additionally, for any $n \in \N$,
\[
\rH^\ast(\rp^n; \Ftwo) \cong \frac{\Ftwo[\sigma]}{(\sigma^{n+1})}.
\]

The action of the Steenrod algebra $\cA_2$ on $\rH^\ast(\rp^n; \Ftwo)$, for $n$ possibly equal to $\infty$, is determined by
\[
\Sq^k(\sigma^\ell) = \binom{\ell}{k}\sigma^{\ell+k},
\]
where the binomial coefficient is evaluated modulo \(2\).

\subsubsection{}\label{sss:steenrod_odd}

For an odd prime \(p\), the Steenrod algebra $\cA_p$ is generated by the \defn{Bockstein} homomorphism of the exact sequence
\[
0 \to \Z/p \to \Z/p^2 \to \Z/p \to 0
\]
and the \defn{Steenrod reduced powers} \(\set{\rP^i}_{i \in \N}\), where $\rP^i \colon \rH^n(\cX; \Fp) \to \rH^{n+2i(p-1)}(\cX; \Fp)$, for any $n\geq 0$.
These are degree \(2i(p-1)\) stable operations which are axiomatically characterized in a similar manner to Steenrod squares.

The \defn{total Steenrod reduced power operation} is the (inhomogeneous) cohomology operation
\[
\rP = \rP^0 + \rP^1 + \rP^2 + \dotsb
\]
acting on \(\bigoplus_{m \in \N} \rH^m(\cX; \Fp)\).

\subsubsection{}\label{sss:cohomology_lens}

Consider the system of unit spheres
\[
\bS^1 \subset \bS^3 \subset \bS^5 \subset \dotsb
\]
defined by the canonical inclusions \(\bC^n \subset \bC^{n+1}\), with the \(\cyc_q\) action, for some \(q \in \N\), defined by the diagonal action of the group of \(q^\th\) roots of unity.
A model for \(K(\Z/q,1)\) is the infinite Lens space \(\rL^\infty_q\), defined as the union of the sequence of (finite) \defn{Lens spaces}
\[
\rL^0_q \subset \rL^1_q \subset \rL^2_q \subset \dotsb,
\]
where each \(\rL^n_q\) is the space of orbits of \(\bS^{2n+1}\) under the \(\cyc_q\)-action.

The mod \(p\) cohomology of \(\rL^\infty_p\) when \(p\) is an odd prime is given by
\[
\rH^\ast(\rL^\infty_p; \Fp)
\cong
\Fp[\sigma] \otimes \Lambda(\tau),
\]
where \(|\sigma| = 2\) and \(|\tau| = 1\).
For finite \(n \in \N\),
\[
\rH^\ast(\rL^n_p; \Fp)
\cong
\frac{\Fp[\sigma]}{(\sigma^{n+1})} \otimes \Lambda(\tau).
\]

The action of the Steenrod algebra \(\cA_p\) on \(\rH^\ast(\rL^n_p; \Fp)\) (for \(n\) possibly equal to \(\infty\)) is determined by:
\[
\beta(\tau) = \sigma,
\qquad
\beta(\sigma) = 0,
\]
and, for all \(\ell \in \N\) and \(k \in \N\),
\[
\rP^k(\sigma^\ell) = \binom{\ell}{k}\sigma^{\ell+k(p-1)},
\qquad
\rP^k(\tau\sigma^\ell) = \binom{\ell}{k}\tau\sigma^{\ell+k(p-1)},
\]
where the binomial coefficient is evaluated modulo \(p\) (and \(\rP^0 = \mathrm{id}\)).


\subsection{Persistent $\theta$-modules}

\subsubsection{}\label{ss:theta-modules}

A cohomology operation $\theta \in \cO(\kk,m; \kk,n)$ is \defn{linear} if \(\k\) is a field and the map \(\theta_{\cX} \colon \rH^m(\cX; \k) \to \rH^n(\cX; \k)\) is a linear transformation for any cellular space \(\cX\).
We denote the subset of $\cO(\kk,m; \kk,n)$ containing the linear cohomology operations by $\cO(m,n;\k)$ and omit the field \(\k\) from the notation when no confusion arises from doing so.
Since every additive map between \(\mathbb F_p\)-vector spaces is automatically \(\mathbb F_p\)-linear, and since the Steenrod operations are additive \cite[\S4.L]{hatcher2000}, every homogeneous operation in \(\cA_p\) of degree \(d\) determines an element of \(\cO(m,m+d;\mathbb F_p)\) for every \(m\).

A linear cohomology operation \(\theta \in \cO(\ell, m)\) induces a morphism of persistence modules, denoted \(\theta_X \colon \rH^\ell(X) \to \rH^m(X)\), for any cellular \(\R\)-space \(X\).
The \defn{image} and \defn{kernel} of $\theta_X$, denoted \defn{$\img_\theta(X)$} and \defn{$\ker_\theta(X)$} respectively, are the persistence modules
\begin{align*}
	\img_\theta(X)_r &= \img((\theta_X)_r)\,, &
	\ker_\theta(X)_r &= \ker((\theta_X)_r)\,, \\
	\img_\theta(X)_{s,t} &= \rH^m(X)_{s,t}\big|_{\img(\theta_t)}\,, &
	\ker_\theta(X)_{s,t} &= \rH^\ell(X)_{s,t}\big|_{\ker(\theta_t)}\,,
\end{align*}
where \(r, s, t \in \R\) and \(s \leq t\).
Here \(\rH^j(X)_{s,t}\) denotes the contravariant map \(\rH^j(X_t) \to \rH^j(X_s)\).

Since the identity and zero maps are linear cohomology operations, for any \(m \in \N\) we have
\[
\img_\id(X) \cong \ker_0(X) \cong \rH^m(X),
\]
which shows that these invariants include persistent cohomology as a special case.

\subsubsection{}\label{ss:theta-modules-q-tame}

If \(\rH^m(X)\) (resp. \(\rH^\ell(X)\)) is q-tame, then \(\img_\theta(X)\) (resp. \(\ker_\theta(X)\)) is also q-tame.
In this case we refer to its barcode as the \defn{$\img_\theta$-barcode} (resp. \defn{\(\ker_\theta\)-barcode}) of \(X\).
The Steenrod barcodes introduced in \cite{medina2022per_st} are specific instances of \(\img_\theta\)-barcodes, where \(\theta\) corresponds to a Steenrod square.


\subsection{Sums and products}\label{ss:sums_products}

We study persistent cohomology operations under sums and products, focusing only on Steenrod operations for the latter.

\subsubsection{}\label{ss:wedge_sum_general}

The \defn{wedge sum} of two pointed \(\R\)-spaces \(X\) and \(Y\) is given by
\[
(X \vee Y)_r = X_r \vee Y_r \quad\text{and}\quad (X \vee Y)_{s,t} = X_{s,t} \vee Y_{s,t},
\]
where \(r, s, t \in \R\) and \(s \leq t\).

\lemma
Let $X$ and $Y$ be pointed cellular $\R$-spaces.
For any linear cohomology operation $\theta$:
\begin{align*}
	\img_\theta(X \vee Y) &\cong \img_\theta(X) \oplus \img_\theta(Y), \\
	\ker_\theta(X \vee Y) &\cong \ker_\theta(X) \oplus \ker_\theta(Y).
\end{align*}

\begin{proof}
	This follows directly from the fact that for any \(r \in \R\), the natural inclusions $X_r \to X_r \vee Y_r \leftarrow Y_r$ induce an isomorphism $\rH^*(X_r \vee Y_r) \to \rH^*(X_r) \oplus \rH^*(Y_r)$ (\cite[Corollary~2.25]{hatcher2000}).
\end{proof}

\subsubsection{}\label{ss:wedge sum}

The \defn{gluing metric} on the wedge sum of pointed metric spaces ${(\cX_i, x_i)}_{i = 0}^k$ is defined as follows \cite{burago2001course}: For $x \in \cX_i$ and $x' \in \cX_j$,
\[
d(x, x') =
\begin{cases}
	d_{\cX_i}(x, x_i) + d_{\cX_j}(x', x_j) & i \neq j, \\
	\hfil d_{\cX_i}(x, x') & i = j.
\end{cases}
\]

\theorem
Let $\cX$ and $\cY$ be pointed metric spaces.
For any linear cohomology operation $\theta$:
\begin{align*}
	\img_\theta^\VR(\cX \vee \cY) &\cong \img_\theta^\VR(\cX) \oplus \img_\theta^\VR(\cY), \\
	\ker_\theta^\VR(\cX \vee \cY) &\cong \ker_\theta^\VR(\cX) \oplus \ker_\theta^\VR(\cY),
\end{align*}
where \(\cX \vee \cY\) is equipped with the gluing metric.

\begin{proof}
	This follows from the previous lemma and the fact, proven in \cite[Proposition~1]{adamaszek2020homotopy}, that the natural inclusion
	\[
	\VR_r(\cX) \vee \VR_r(\cY) \to \VR_r(\cX \vee \cY)
	\]
	is a homotopy equivalence for each $r \in \R$ when \(\cX \vee \cY\) is equipped with the gluing metric.
\end{proof}

\subsubsection{}\label{ss:product_general}

Similarly, for any two $\R$-spaces $X$ and $Y$, the $\R$-space $X \times Y$, their \defn{product}, is defined by
\[
(X \times Y)_r = X_r \times Y_r \quad\text{and}\quad (X \times Y)_{s,t} = (X_{s,t} \times Y_{s,t}),
\]
where \(r, s, t \in \R\) and \(s \leq t\).
Analogously, the \defn{tensor product} \(V \ot W\) of two persistence modules is defined by
\[
(M \ot N)_r = M_r \ot N_r \quad\text{and}\quad (M \ot N)_{s,t} = (M_{s,t} \ot N_{s,t}),
\]
where \(r, s, t \in \R\) and \(s \leq t\).

\lemma
Let $X$ and $Y$ be pointed cellular $\R$-spaces, at least one of them landing in the full subcategory of finite type spaces.
Over any prime \(p\):
\begin{align*}
	\img_\rP(X \times Y) &\cong \img_\rP(X) \ot \img_\rP(Y),\\
	\ker_\rP(X \times Y) &\cong \ker_\rP(X) \ot \ker_\rP(Y),
\end{align*}
where \(\rP\) is the total Steenrod power if \(p > 2\) and \(\rP = \Sq\) if \(p = 2\), acting on unreduced cohomology.

\begin{proof}
	For any \(r \in \R\), the natural projections $X_r \leftarrow X_r \times Y_r \to Y_r$ induce the K\"unneth isomorphism $\rH^*(X_r) \otimes \rH^*(Y_r) \to \rH^*(X_r \times Y_r)$, since we are working with cellular spaces over a field. (\cite[Theorem~3.15]{hatcher2000}).
	Under this isomorphism, the (exterior) Cartan formula identifies \(\rP_{X_r \times Y_r}\) with \(\rP_{X_r} \otimes \rP_{Y_r}\) (\cite[Corollary~2.7]{may1970general}).
	Over a field, the image of a tensor product of linear maps is the tensor product of their images, proving the image formula.

	For the kernel formula, note that \(\rP\) is injective: since \(\rP^0 = \id\) and every other component raises degree.
	Thus, all three kernels in the formula are zero.
	The resulting identifications are natural in \(r\), so they define isomorphisms of persistence modules.
\end{proof}

\subsubsection{}\label{ss:products}

Given two metric spaces $\cX$ and $\cY$.
The \defn{$\ell^\infty$ metric} on $\cX \times \cY$ is defined by
\[
d((x,y), \, (x',y')) = \max\set{d_{\cX}(x,x'), d_{\cY}(y,y')}.
\]

In \cite[Proposition~10.2]{adamaszek2017VietorisProduct}, it was shown that the Vietoris--Rips complex of a product space is homotopy equivalent to the product of the Vietoris--Rips complexes of the individual spaces at any specific scale.
Furthermore, the proof of \cite[Theorem~6.1(1)]{lim2024vietoris} constructs a natural isomorphism $\rH^\VR(\cX) \ot \rH^\VR(\cY) \cong \rH^\VR(\cX \times \cY)$, which, combined with the previous lemma gives the following.

\theorem
Let $\cX$ and $\cY$ be metric spaces.
Over any prime:
\begin{align*}
	&\img^{\VR}_\rP(\cX \times \cY) \cong \img^{\VR}_\rP(\cX) \ot \img^{\VR}_\rP(\cY), \\
	&\ker^{\VR}_\rP(\cX \times \cY) \cong \ker^{\VR}_\rP(\cX) \ot \ker^{\VR}_\rP(\cY),
\end{align*}
where \(\rP\) is the total Steenrod power if \(p > 2\), \(\rP = \Sq\) if \(p = 2\), acting on unreduced cohomology, and \(\cX \times \cY\) is equipped with the \(\ell^\infty\) metric.

\subsection{Stability}\label{ss:stability}

We prove the stability of cohomology operations using the homotopy interleaving distance of \cite{blumberg2023interleaving}.

\subsubsection{}\label{lem:w.h.e. preservance}

\lemma If two cellular $\R$-spaces $X$ and $X'$ are weakly equivalent (\cref{ss:R-spaces}), then, for any linear cohomology operation $\theta$
\[
\begin{split}
	\img_\theta(X) & \cong \img_\theta(X'), \\
	\ker_\theta(X)& \cong \ker_\theta(X').
\end{split}
\]

\begin{proof}
	We write the proof only for $\img_\theta(X)$ since that for $\ker_\theta(X)$ is analogous.
	Given that $X$ and $X'$ are weakly equivalent, there exists an $\R$-space $Z$ and morphisms $f \colon Z \to X$ and $g \colon Z \to X'$ such that $f_r$ and $g_r$ are weak homotopy equivalences for any $r \in \R$.
	Without loss of generality we can assume $Z$ to be cellular using the cofibrant replacement functor of the model category structure on $\R$-spaces described in \cite{blumberg2023interleaving}.
	Both $f_r$ and $g_r$ induce an isomorphism on cohomology for each $r \in \R$, and the naturality of $\theta$ implies that $\img_\theta(Z)$ and $\img_\theta(X)$ are isomorphic, as well as $\img_\theta(Z)$ and $\img_\theta(X')$.
	This finishes the proof.
\end{proof}

\subsubsection{}\label{thm:stability theta}

\lemma
For cellular $\R$-spaces $X$ and $Y$ and a linear cohomology operation $\theta$
\begin{align*}
	\di(\img_\theta(X), \img_\theta(Y)) &\leq \dhi(X,Y), \\
	\di(\ker_\theta(X), \ker_\theta(Y)) &\leq \dhi(X,Y).
\end{align*}

\begin{proof}
	We only write the proof for $\img_\theta(X)$ since that for $\ker_\theta(X)$ is similar.
	Take any $\delta > \dhi(X,Y)$.
	By the definition of the homotopy interleaving distance, there exist $\R$-spaces $X' \simeq X$ and $Y' \simeq Y$ such that $X'$ and $Y'$ are $\delta$-interleaved.
	Applying the cofibrant replacement functor $\rQ$ we get $d_\rI(\rQ X', \rQ Y') \leq \delta$ since applying any functor does not increase the interleaving distance.
	Both $\rQ X$ and $\rQ Y$ can be assumed to be cellular $\R$-spaces, since cofibrant objects in $\Top^\R$ are in particular objectwise cofibrant.
	By applying the triangle inequality and Lemma~\ref{lem:w.h.e. preservance}, we obtain
	\begin{align*}
		\di(\img_\theta(X), \img_\theta(Y)) \leq& \,
		\di(\img_\theta(X), \img_\theta(\rQ X')) \\ +& \,
		\di(\img_\theta(\rQ X'), \img_\theta(\rQ Y')) + \di(\img_\theta(\rQ Y'), \img_\theta(Y)) \\ =& \,
		0 + \di(\img_\theta(\rQ X'), \img_\theta(\rQ Y')) + 0 \\ \leq \,&
		\di(\rQ X', \rQ Y') \\ \leq \,&
		\delta.
	\end{align*}
	where, for the second inequality, we used Proposition~\ref{ss:interleaving}.
	Since $\delta > \dhi(X,Y)$ is arbitrary, we obtain the desired inequality.
\end{proof}

\subsubsection{}\label{thm:stability VR}

For a metric space \(\cX\) and a linear cohomology operation \(\theta\) we denote \(\img_\theta\big(\VR(\cX)\big)\) and \(\ker_\theta\big(\VR(\cX)\big)\) respectively by \(\img_\theta^\VR(\cX)\) and \(\ker_\theta^\VR(\cX)\).

\theorem
For any two metric spaces $\cX$ and $\cY$ and a linear cohomology operation~$\theta$
\begin{align*}
	\di(\img_\theta^\VR(\cX),\, \img_\theta^\VR(\cY)) \leq 2 \cdot \dgh(\cX,\cY), \\
	\di(\ker_\theta^\VR(\cX),\, \ker_\theta^\VR(\cY)) \leq 2 \cdot \dgh(\cX,\cY).
\end{align*}

\begin{proof}
	From \cref{thm:stability-HI} and the previous theorem we obtain
	\[
	\di\big(\img_\theta^\VR(\cX), \img_\theta^\VR(\cY)\big) \leq
	\dhi\big(\VR(\cX), \VR(\cY)\big) \leq 2 \cdot \dgh(\cX,\cY).
	\]
	The same argument applies to $\ker_\theta^\VR$.
\end{proof}

We remark that if the metric spaces are compact or, more generally, totally bounded, the above persistent modules are q-tame and the bottleneck distance of their barcodes is bounded by twice the Gromov--Hausdorff distance, since bottleneck and interleaving distances agree.

\section{Critical radii and spherical quotients}\label{s:barcodes}

In this section, we utilize Gromov's concept of the filling radius and its modifications to provide estimates for the Vietoris--Rips barcodes and \(\theta\)-barcodes of Riemannian pseudomanifolds.
The primary examples analyzed include wedge sums and quotients of round spheres.


\subsection{Kuratowski embedding and critical radii}\label{sub:filling radii}

\subsubsection{}

Any compact metric space \(\cX\) can be isometrically embedded into the space \(\rL^\infty(\cX)\) of all bounded real-valued functions on \(\cX\) via the map \(x \in \cX \mapsto d_\cX(x,\cdot)\), known as the \defn{Kuratowski embedding}.
The \(\R\)-space obtained by considering radius \(r\) neighborhoods of the image of \(\cX\) is denoted by \(\rU_r(\cX)\).

Throughout the text, Riemannian manifolds are thought of as metric spaces with the geodesic distance.

\subsubsection{}\label{ss:filling_radii_def}

Following \cite{gromov1983filling}, the \defn{filling radius} of an \(n\)-dimensional closed Riemannian manifold \(\cM\), denoted \(\fillrad(\cM)\), is defined as the infimal \(\epsilon > 0\) such that the fundamental class in \(\rH_n(\cM; G)\) vanishes under the inclusion \(\cM \hookrightarrow \rU_\epsilon(\cM)\), where \(\rU_\epsilon(\cM)\) is the \(\epsilon\)-neighborhood of \(\cM\) in \(\rL^\infty(\cM)\) and $G$ is an abelian group.

In this work, we consider the filling radius of \(\cM\) with coefficients in a field \(\field\) as
\[
\fillrad(\cM;\field)
=
\min \big\{ r \mid (0,r) \in \barc \rH_n(\rU(\cM);\field) \big\}.
\]
For non-orientable manifolds (eg, \(\mathbb{RP}^n\)), we take \(\field = \Ftwo\) so that the top-dimensional fundamental class is nontrivial.
For spheres $S^n$, the choice of coefficients does not affect the filling radius.
For lens spaces $\rL_p^n$ (with $p$ an odd prime), while any field suffices to define the filling radius, we work over $\field = \Fp$ to detect the $p$-torsion in degree one.
We adopt these conventions throughout and omit $\field$ from the notation whenever it is clear from the context.

\subsubsection{}\label{ss:first_critical_value}

Let \(\cM\) be a closed Riemannian manifold.
By \cite[Theorem~3.5]{hausmann1995vietoris} and \cite[Theorem~4.1]{lim2024vietoris}, we know that there exists \(r_\cM > 0\) such that for all \(r \in (0, r_\cM)\) the inclusion \(\cM \to \rU_{r}(\cM)\) is a homotopy equivalence.
This leads to considering the supremum of all such \(r > 0\) for which the inclusion \(\cM \to \rU_r(\cM)\) is a homotopy
equivalence.
We will refer to this value as the homotopical radius of $\cM$ and denote it by \(\crit(\mathcal{M})\).

\subsubsection{} \label{ss:beta v.s. fillrad}

In \cite[Definition~9.44]{lim2024vietoris}, the authors define the generalized filling radius of a homology class \(\omega\) in \(X\) as the smallest \(r > 0\) for which \(\omega\) becomes trivial in \(\rU_r(\cX)\), capturing the scale at which \(\omega\) vanishes under the inclusion of \(X\) into its \(r\)-neighborhood.
Below, we consider a weaker version of this concept, which is equivalent to the smallest generalized filling radius among all non-zero homology classes \(\omega\) in \(X\) of a given degree $m$.

Let us fix a field \(\field\) and an integer \(m \in \N\), the \defn{\(m\)-homological radius} of a totally bounded metric space \(\cX\) is defined by
\[
\firstdeath{m}{\cX; \field} =
\min \big\{r \mid (0, r) \in \barc \rH_m(\rU(\cX); \field)\big\},
\]
if \(\rH_m(\cX) \neq 0\).
If \(\rH_m(\cX) = 0\), set \(\firstdeath{m}{\cX} = 0\).
The \(m\)-homological radius may depend on the coefficient field.
Throughout, we follow the convention in \cref{ss:filling_radii_def}, and omit \(\field\) from the notation whenever it is clear from context.

Similarly, for a linear cohomology operation \(\theta\), we define the \defn{\(\img_\theta\)-radius} of \(\cX\) by
\[
\firstdeath{\theta}{\cX} = \min \big\{r \mid (0, r) \in \barc \img_\theta(\rU(\cX))\big\},
\]
if \(\img_\theta(\cX) \neq 0\).
If \(\img_\theta(\cX) = 0\), set $\firstdeath{\theta}{\cX} = 0$.
The \(\ker_\theta\)-radius is defined similarly, but we do not use it in this work.

Clearly, if \(\rH_m(\cX) \neq 0\) (resp. \(\img_\theta(\cX) \neq 0\)), then
\[
\crit(\cX) \leq \firstdeath{m}{\cX} \qquad (\text{resp. }\crit(\cX) \leq \firstdeath{\theta}{\cX}),
\]
and if $\cM$ is connected and $n$-dimensional, then
\[
\firstdeath{n}{\cM} = \fillrad(\cM).
\]

\subsubsection{}\label{ss:kuratowski_vr}

There is a deep relationship between the \(r\)-neighborhood filtration of $\cX$ in $\rL^\infty(\cX)$ and the Vietoris--Rips complex of \(\cX\) at scale \(2r\), stated in the following.

\proposition \cite[Theorem~4.1]{lim2024vietoris} The spaces $\VR_{2r}(\cX)$ and $\rU_r(\cX)$ are naturally homotopy equivalent for every \(r > 0\).

\medskip Due to the natural homotopy equivalence between these filtrations, we will use them interchangeably as needed.
In particular, this equivalence allows us to bound the Vietoris--Rips barcode of Riemannian pseudomanifolds using their critical radii.

\subsection{General estimates}\label{ss:barcode_general}

Let \(\cM\) be a closed Riemannian manifold.
Consider \(\ell, m \in \N\) and \(\theta \in \cO(\ell,m)\).
We will simplify notation writing \(R = \diam(\cM),\, \alpha = 2\crit(\cM)\), \(\beta_m = 2\firstdeath{m}{\cM}\), and \(\gamma_\theta = 2\firstdeath{\theta}{\cM}\).
We say a bar $(a, b)$ is \defn{dominated} by another $(c,d)$ if $c \leq a < b \leq d$.

Because $\VR_r(\cM)$ is empty for \(r \leq 0\) and the homotopy type of $\VR_r(\cM)\simeq \rU_{r/2}(\cM)$ remains the homotopy type of $\cM$ for $r \in (0, \alpha)$, bars in \(\barc \rH_m^{\VR}(\cM)\) and $\barc\img_\theta^{\VR}(\cM)$ either start at $0$ or start after $\alpha$.

More precisely, if \(k = \dim \opH_m(\cM) > 0\) then $\barc\rH_m^{\VR}(\cM)$ contains $(0, \beta_m)$ and \((k - 1)\) additional bars of the form \((0, b)\) with \(\beta_m \leq b \leq R\).
Additionally, all other bars are dominated by \((\alpha, R)\).
If \(\dim \opH_m(\cM) = 0\) then all bars in \(\barc\rH_m^{\VR}(\cM)\) are dominated by \((\alpha, R)\).
See the first row of \cref{fig:barcodes_general} for a pictorial representation of these estimates.

A similar analysis applies to $\barc\img_\theta^{\VR}(\cM)$.
See the second row of \cref{fig:barcodes_general} for these estimates.

\begin{figure}
	\centering
	\begin{tikzpicture}[scale=0.52]
	\begin{axis} [
		title = {\LARGE $\barc\opH_m^{\VR}(\cM)$, if $\opH_m(\cM) \neq 0$},
		ticklabel style = {font=\Large},
		axis y line=middle,
		axis x line=middle,
		ytick={0.7,0.95},
		yticklabels={$2\firstdeath{m}{\cM}$,$\diam(\cM)$},
		xtick={0.55,0.95},
		xticklabels={$2\crit(\cM)$, $\diam(\cM)$},
		xmin=-0.015, xmax=1.1,
		ymin=0, ymax=1.1,]
		\addplot [mark=none] coordinates {(0,0) (1,1)};
		\addplot [thick,color=black!20!white,fill=black!30!white,
		fill opacity=0.4]coordinates {
			(0.55,0.95)
			(0.55,0.55)
			(0.95,0.95)
			(0.55,0.95)};
		\addplot [black!40!white,mark=none,dashed, thin] coordinates {(0,0.7) (0.7,0.7)};
		\addplot [black!40!white,mark=none,dashed, thin] coordinates {(0,0.55) (0.55,0.55)};
		\addplot [black!40!white,mark=none,dashed, thin] coordinates {(0.55,0) (0.55,0.55)};
		\addplot[line width=1.5mm, color=black!30!white] coordinates{(0, 0.7) (0, 0.95)};
		\addplot[barccolor,mark=*] (0, 0.7) circle (2pt) node[above right,barccolor]{};
	\end{axis}
\end{tikzpicture}
\begin{tikzpicture}[scale=0.52]
	\begin{axis} [
		title = {\LARGE $\barc\opH_m^{\VR}(\cM)$, if $\opH_m(\cM) = 0$},
		ticklabel style = {font=\Large},
		axis y line=middle,
		axis x line=middle,
		ytick={0.95},
		yticklabels={$\diam(\cM)$},
		xtick={0.55,0.95},
		xticklabels={$2\crit(\cM)$, $\diam(\cM)$},
		xmin=-0.015, xmax=1.1,
		ymin=0, ymax=1.1,]
		\addplot [mark=none] coordinates {(0,0) (1,1)};
		\addplot [thick,color=black!20!white,fill=black!30!white,
		fill opacity=0.4]coordinates {
			(0.55,0.95)
			(0.55,0.55)
			(0.95,0.95)
			(0.55,0.95)};
		\addplot [black!40!white,mark=none,dashed, thin] coordinates {(0,0.55) (0.55,0.55)};
		\addplot [black!40!white,mark=none,dashed, thin] coordinates {(0.55,0) (0.55,0.55)};
	\end{axis}
\end{tikzpicture}

\begin{tikzpicture}[scale=0.52]
	\begin{axis} [
		title={\LARGE $\thetabarc{\cM}$, if $\img\theta_{\cM}\neq 0$},
		ticklabel style = {font=\Large},
		axis y line=middle,
		axis x line=middle,
		ytick={0.7,0.95},
		yticklabels={$2\firstdeath{\theta}{\cM}$,$\diam(\cM)$},
		xtick={0.55,0.95},
		xticklabels={$2\crit(\cM)$, $\diam(\cM)$},
		xmin=-0.015, xmax=1.1,
		ymin=0, ymax=1.1,]
		\addplot [mark=none] coordinates {(0,0) (1,1)};
		\addplot [thick,color=black!20!white,fill=black!30!white,
		fill opacity=0.4]coordinates {
			(0.55,0.95)
			(0.55,0.55)
			(0.95,0.95)
			(0.55,0.95)};
		\addplot [black!40!white,mark=none,dashed, thin] coordinates {(0,0.7) (0.7,0.7)};
		\addplot [black!40!white,mark=none,dashed, thin] coordinates {(0,0.55) (0.55,0.55)};
		\addplot [black!40!white,mark=none,dashed, thin] coordinates {(0.55,0) (0.55,0.55)};
		\addplot[line width=1.5mm, color=black!30!white] coordinates{(0, 0.7) (0, 0.95)};
		\addplot[barccolor,mark=*] (0, 0.7) circle (2pt) node[above right,barccolor]{};
	\end{axis}
\end{tikzpicture}
\begin{tikzpicture}[scale=0.52]
	\begin{axis} [
		title={\LARGE $\thetabarc{\cM}$, if $\img\theta_{\cM}= 0$},
		ticklabel style = {font=\Large},
		axis y line=middle,
		axis x line=middle,
		ytick={0.95},
		yticklabels={$\diam(\cM)$},
		xtick={0.55,0.95},
		xticklabels={$2\crit(\cM)$, $\diam(\cM)$},
		xmin=-0.015, xmax=1.1,
		ymin=0, ymax=1.1,]
		\addplot [mark=none] coordinates {(0,0) (1,1)};
		\addplot [thick,color=black!20!white,fill=black!30!white,
		fill opacity=0.4]coordinates {
			(0.55,0.95)
			(0.55,0.55)
			(0.95,0.95)
			(0.55,0.95)};
		\addplot [black!40!white,mark=none,dashed, thin] coordinates {(0,0.55) (0.55,0.55)};
		\addplot [black!40!white,mark=none,dashed, thin] coordinates {(0.55,0) (0.55,0.55)};
	\end{axis}
\end{tikzpicture}
	\caption{Let $\cM$ be a closed Riemannian manifold.
		\emph{Top row:} persistent reduced homology barcodes of $\cM$.
		\emph{Bottom row:} $\img_\theta$-barcodes of $\cM$.
		In each figure, the gray region represents where additional bars could potentially exist within the corresponding barcode.
		See \cref{ss:barcode_general} for details.}
	\label{fig:barcodes_general}
\end{figure}


\subsection{Spheres and their wedge sum}\label{ss:Sn}

For any integer $n \geq 1$ and real number $\rho > 0$, let $\bS^n(\rho)$ be the \defn{$n$-sphere of radius $\rho$}
\[
\bS^n(\rho) = \set{x \in \R^{n+1} \mid \bars{x} = \rho}.
\]
With the induced Riemannian metric, \(\bS^n(\rho)\) is referred to as a \defn{round sphere}.
We simplify notation writing \(\bS^n\) instead of \(\bS^n(1)\).
The constant
\[
\zeta_n = \arccos\Big(\frac{-1}{n+1}\Big),
\]
which is the diameter of an inscribed regular $(n+1)$-simplex in $\bS^n$, will play an important role.

\subsubsection{}\label{ss:critical values of Sn}

Characterizing the homotopy types of Vietoris--Rips complexes is a challenging question, even for simple spaces like the $n$-spheres.
For the circle $\bS^1$, a complete answer is provided in \cite{adamaszek2017VietorisProduct}, but for more general spheres $\bS^n$ with $n>1$, only partial results are known.

In \cite{adamaszek2018metric}, the authors introduced the Vietoris--Rips metric thickening, a metrization of the Vietoris--Rips complex that embeds the original metric space as a metric subspace. They determined the first scale parameter at which the homotopy type of the thickening of an \(n\)-sphere changes.
While this does not directly yield the homotopical radius of the \(n\)-sphere, it provides insights into its homotopy behavior and yields homological information, since the thickening shares the same homology as the original Vietoris–Rips complex \cite[Corollary~3]{gillespie2024vietoris}.

Building on the equivalence between the Kuratowski and Vietoris--Rips filtrations, the authors in \cite{lim2024vietoris} determined both the homotopical radius and the subsequent homotopy type in the Vietoris--Rips filtration of the \(n\)-sphere.

The following lemma is a consequence of these results, along with those from \cite{katz1983filling}.

\lemma
For any $m,n \in \N$ and linear cohomology operation $\theta$ of non-zero degree we have:
\begin{enumerate}
	\item \(\crit(\bS^n) = \frac{1}{2}\zeta_n\),
	\item \(\firstdeath{m}{\bS^n} =
	\begin{cases}
		\frac{1}{2}\zeta_n & m = n, \\
		\hfil 0 & m \neq n,
	\end{cases}\)
	\item \(\firstdeath{\theta}{\bS^n} = 0\).
\end{enumerate}

\begin{proof}
	(1) By \cite[Theorem~7.1]{lim2024vietoris}, for any $0 < r \leq \zeta_n$, the space $\VR_r(\bS^n)$ is homotopy equivalent to $\bS^n$, and the homotopy type of $\VR_r(\bS^n)$ changes at $\zeta_n$.
	This implies $\crit(\bS^n)=\frac{1}{2}\zeta_n$.

	\smallskip (2) According to \cite{katz1983filling}, \(\fillrad(\bS^n) = \frac{1}{2}\zeta_n\).
	Applying \cref{ss:beta v.s. fillrad}, we obtain \(\firstdeath{n}{\bS^n} = \frac{1}{2}\zeta_n\).
	When $m\neq n$, the statement follows from the fact that the initial space in the filtration has a trivial degree $m$ reduced homology.

	\smallskip (3) We apply a similar argument as in the $m\neq n$ case of (2). The statement follows from the fact that the initial space in the filtration has trivial $\img_\theta$.
\end{proof}

\subsubsection{}\label{ss:VRSn projection}

Let \(n \in \N\) and \(r \in (0, \zeta_n]\).
From \cite[Theorem~7.1]{lim2024vietoris} we know that the spaces \(\VR_r(\bS^n)\) and \(\bS^n\) have the same homotopy type for any \(n \in \N\).
We now recall from \cite{adamaszek2018metric} a natural map that can be used to realize this equivalence.

In \cite{adamaszek2018metric}, the authors introduced the concept of Vietoris–Rips metric thickenings, which endows the Vietoris–Rips complex with a new topology.
They define the \defn{canonical projection of \(\VR_r(\bS^n)\)} as the set map
\[
f_r^n \colon \VR_r(\bS^n) \to \R^{n+1} \setminus \set{0} \to \bS^n,
\]
constructed by first mapping a formal linear combination \(\sum \lambda_i x_i\) to the corresponding point \(\sum \lambda_i x_i\) in \(\mathbb{R}^{n+1}\), followed by a radial projection onto \(\bS^n\).
Throughout this work, we do not distinguish between a simplicial complex and its geometric realization.

As mentioned in \cite{adamaszek2018metric}, the map $f_r^n$ is well-defined because, when $r \leq \zeta_n$, the map $\VR_r(\bS^n) \to \R^{n+1}$ misses the origin by the proof of \cite[Lemma 3]{lovasz1983self}.
Moreover, the authors proved that when \(\VR_r(\bS^n)\) is equipped with the metric thickening topology, the canonical projection is a homotopy equivalence.

More recently, \cite{gillespie2024vietoris} showed that the identity map between the Vietoris--Rips complex and the Vietoris--Rips metric thickening is a weak equivalence.
Consequently, the canonical projection of \(\VR_r(\bS^n)\) is also a weak equivalence for the usual Vietoris--Rips complex. By the Whitehead theorem, this weak equivalence between CW spaces is indeed a homotopy equivalence.
Thus, the canonical projection remains a homotopy equivalence when \(\VR_r(\bS^n)\) is equipped with the usual topology, as summarized by the following lemma.

\lemma
The canonical projection of the Vietoris--Rips complex \(\VR_r(\bS^n)\) is a homotopy equivalence for every $r \in (0, \zeta_n]$.

\subsubsection{}\label{ss:barcode_Sn}

We apply \cref{ss:critical values of Sn} and the barcode estimates from \cref{ss:barcode_general} to \(\bS^n\).
Using the facts that (1) the reduced homology of \(\bS^n\) has dimension one for \(m = n\) and is zero for all other values of \(m\), and (2) for any linear cohomology operation \(\theta \in \cO(\ell, m)\) with \(\ell \neq m\), \(\img \theta_{\bS^n}\) is trivial, we proceed with the estimation but omit the figure, as it is simply a special case of \cref{fig:barcodes_vs},
corresponding to $u_n=1$ and $u_i=0$ for $i\neq n$.

\subsubsection{}\label{ss:barcodes_VS}

For $n \in \N$ and $u_1, \dots, u_n \in \N^n$, let
\[
\VS^{u_1,\dots,u_n} =
\overbrace{\bS^1\vee\dots\vee\bS^1}^{u_1} \vee\dots\vee \overbrace{\bS^n\vee\dots\vee\bS^n}^{u_n}.
\]
Using the homotopy equivalence between the Vietoris--Rips complex of a metric gluing with the wedge sum of the Vietoris--Rips complex described in \cref{ss:wedge sum}, we have the following isomorphisms of persistence modules for \(m \in \N\) and \(\theta \in \cO(\ell,m)\):
\[
\begin{split}
	\rH_m^\VR(\VS^{u_1,\dots,u_n}) &\cong \, \bigoplus_{i=1}^n \rH_m^\VR(\bS^i)^{\oplus u_i}, \\
	\img_\theta^\VR(\VS^{u_1,\dots,u_n}) &\cong \, \bigoplus_{i=1}^n \img_\theta^\VR(\bS^i)^{\oplus u_i}.
\end{split}
\]

Therefore, both the homology barcodes and \(\img_\theta\)-barcodes of \(\VR(\VS^{u_1, \dots, u_n})\) are the multiset unions of the corresponding barcodes of \(\bS^i\); see \cref{fig:barcodes_vs}.
In particular, each $m$-sphere summand contributes a bar $(0,\zeta_m)$
to $\Hbarc[\field]{m}{\VS^{u_1,\dots,u_n}}$,
and this occurs with multiplicity $u_m$.
Since $\zeta_1 > \dots > \zeta_n$,
the smallest threshold $\zeta_n$ determines the vertical dashed line
in Figure~\ref{fig:barcodes_vs}.

\begin{figure}
	\centering
	\begin{tikzpicture}[scale=0.52]
	\begin{axis} [
		title = {\LARGE $\Hbarc[\field]{m}{\VS},\, u_m\neq 0$},
		ticklabel style = {font=\Large},
		axis y line=middle,
		axis x line=middle,
		ytick={0.5,0.6,0.67,0.95},
		yticklabels={,$ \zeta_m$,,$\pi$},
		xtick={0.5,0.55,0.95},
		xticklabels={ ,$\zeta_n$, $\pi$},
		xmin=-0.015, xmax=1.1,
		ymin=0, ymax=1.1,]
		\addplot [mark=none] coordinates {(0,0) (1,1)};
		\addplot [thick,color=black!20!white,fill=black!30!white,
		fill opacity=0.4]coordinates {
			(0.55,0.95)
			(0.55,0.55)
			(0.95,0.95)
			(0.55,0.95)};
		\addplot [black!40!white,mark=none,dashed, thin] coordinates {(0,0.6) (0.6,0.6)};
		\addplot [black!40!white,mark=none,dashed, thin] coordinates {(0,0.55) (0.55,0.55)};
		\addplot [black!40!white,mark=none,dashed, thin] coordinates {(0.55,0) (0.55,0.55)};
		\addplot[barccolor,mark=*] (0, 0.6) circle (2pt) node[above right,barccolor]{};
	\end{axis}
\end{tikzpicture}
\begin{tikzpicture}[scale=0.52]
	\begin{axis} [
		title = {\LARGE $\Hbarc[\field]{m}{\VS},\, m>n$ or $u_m = 0$},
		ticklabel style = {font=\Large},
		axis y line=middle,
		axis x line=middle,
		ytick={0.5,0.55,0.67,0.95},
		yticklabels={,$\zeta_n$,,$\pi$},
		xtick={0.5,0.55,0.95},
		xticklabels={ ,$\zeta_n$, $\pi$},
		xmin=-0.015, xmax=1.1,
		ymin=0, ymax=1.1,]
		\addplot [mark=none] coordinates {(0,0) (1,1)};
		\addplot [thick,color=black!20!white,fill=black!30!white,
		fill opacity=0.4]coordinates {
			(0.55,0.95)
			(0.55,0.55)
			(0.95,0.95)
			(0.55,0.95)};
		\addplot [black!40!white,mark=none,dashed, thin] coordinates {(0,0.55) (0.55,0.55)};
		\addplot [black!40!white,mark=none,dashed, thin] coordinates {(0.55,0) (0.55,0.55)};
	\end{axis}
\end{tikzpicture}

\begin{tikzpicture}[scale=0.52]
	\begin{axis} [
		title={\LARGE $\thetabarc{\VS},\, \theta \in \cO(\ell,m)$},
		ticklabel style = {font=\Large},
		axis y line=middle,
		axis x line=middle,
		ytick={0.5,0.55,0.67,0.95},
		yticklabels={,$\zeta_n$,,$\pi$},
		xtick={0.5,0.55,0.95},
		xticklabels={ ,$\zeta_n$, $\pi$},
		xmin=-0.015, xmax=1.1,
		ymin=0, ymax=1.1,]
		\addplot [mark=none] coordinates {(0,0) (1,1)};
		\addplot [thick,color=black!20!white,fill=black!30!white,
		fill opacity=0.4]coordinates {
			(0.55,0.95)
			(0.55,0.55)
			(0.95,0.95)
			(0.55,0.95)};
		\addplot [black!40!white,mark=none,dashed, thin] coordinates {(0,0.55) (0.55,0.55)};
		\addplot [black!40!white,mark=none,dashed, thin] coordinates {(0.55,0) (0.55,0.55)};
	\end{axis}
\end{tikzpicture}
	\caption{Let $\VS = \VS^{u_1,\dots,u_n}$ for some tuple of non-negative integers.
		\emph{Top row:} persistent reduced homology barcodes of $\VS$, where the dot $(0,\zeta_m)$ has multiplicity $u_m$ which can be zero.
		\emph{Bottom row:} $\img_\theta$-barcodes of $\VS$ where $\theta \in \cO(\ell,m)$ with \(\ell \neq m\).
		In each figure, the gray region represents where additional bars could exist within the corresponding barcode.
		See \cref{ss:barcode_Sn} for the case when the wedge sum is a single sphere and see \cref{ss:barcodes_VS} for the general case.}
	\label{fig:barcodes_vs}
\end{figure}


\subsection{General quotients}

In this section, we review results from \cite{adams2022metric} concerning specific group actions on metric spaces, emphasizing how these actions interact well with the Vietoris--Rips filtration.
Although the original paper primarily focuses on Vietoris--Rips thickenings, we leverage their findings to study the Vietoris--Rips filtration of metric spaces with group actions in later sections.

\subsubsection{}

Let $G$ be a group acting on a metric space $\cX$.
We denote its orbit space, equipped with the quotient topology, by $\cX_G$.
For $x \in \cX$, its orbit will be denoted by $[x]$.
The action of $G$ is called \defn{proper} if, for every $x \in \cX$, there is some $r > 0$ such that $\{g \mid g\cdot B(x,r) \cap B(x,r) \neq \emptyset\}$ is finite.
We recall that $G$ \defn{acts by isometries} on $\cX$, if the map $g \colon \cX \to \cX$ is an isometry for every $g \in G$.

Let $G$ be a group acting properly and by isometries on a metric space $\cX$.
Then the \defn{quotient metric}, defined by
\[
d_{\cX_G}\big([x], [x']\big) = \inf_{g \in G} d_\cX(x, g \cdot x'),
\]
is well-defined on $\cX_G$.

\subsubsection{}\label{ss:h}

Let $G$ be a group acting properly and by isometries on a metric space $\cX$.
The \(G\)-action on $\cX$ induces a natural \(G\)-action on the Vietoris--Rips complex $\VR(\cX)$.
Explicitly, for any \(r > 0, g\in G\) and $\sum \lambda_i x_i \in \VR_r(\cX)$,
\[
g \cdot \sum \lambda_i x_i = \sum \lambda_i (g \cdot x_i).
\]
With respect to this action, there is an induced map on orbit spaces given by
\begin{align*}
	h_r \colon \VR_r(\cX)_G &\to  \VR_r(\cX_G) \\
	 \textstyle \big[\sum\lambda_i x_i\big] & \mapsto \textstyle\sum\lambda_i [x_i].
\end{align*}
If the group action satisfies the \textit{strong \(r\)-diameter} condition, a concept we review below, then for any \(s \leq r\), the map $h_s$ is an isomorphism of simplicial complexes.
This result was originally stated in \cite[Proposition~3.5]{adams2022metric} under an assumption weaker than strong \(r\)-diameter.
We learned about a counterexample to the original statement and a stronger condition that addresses it, identified by Liam Barham, through discussions with the authors of \cite{adams2022metric}, who informed us of Barham's forthcoming work, subsequently presented in \cite{barham2025group}, and facilitated our communication with him.

\subsubsection{}\label{sss:strong_r_action}

The action of a group $G$ on $\cX$ is said to be a \defn{strong $r$-diameter action}, where $r > 0$, if for any non-negative integer $k$, the condition $\diam_{\cX_G}\set[\big]{[x_0],\dots,[x_k]} < r$ implies the existence of a unique choice of $g_i$ for each $i \in \set{1,\dots,k}$ such that $\diam_{\cX}\set{x_0, g_1x_1, \dots, g_kx_k} < r$.
It is not hard to see that a strong $r$-diameter action is free.

Recall from \cref{sss:cohomology_rpn} and \cref{sss:cohomology_lens} the canonical $\cyc_2$- and $\cyc_q$-actions on spheres, for $q \in \N$, that define real projective spaces and Lens spaces, respectively.

The case of $\cyc_2$ in the following lemma was established in \cite[Corollary~4.3]{adams2022metric}.
We extend the proof to $\cyc_q$ for \(q \geq 2\).

\lemma
The canonical $\cyc_2$ and $\cyc_q$ actions on round unit spheres are strong $r$-diameter for every $0 < r \leq \tfrac{2\pi}{q(q+1)}$.

\begin{proof}
	Let $\bS$ be a round unit sphere of odd-dimension (or of any dimension if $q=2$), $\rL_q$ be the Lens space defined by the action of \(\cyc_q\) on \(\bS\), and $\omega = e^{2\pi i/q}$.
	Assume $\diam_{\rL_q}\{[x_0], \dots, [x_k]\} < r$.
	Then, for each $i \in \{1, \dots, k\}$, there exists $g_i \in \cyc_q$ such that $d_{\bS}(x_0, g_i x_i) < r$. For each $i$, let $x_i^* = g_i x_i$, and thus $d_{\bS}(x_0, x_i^*) < r$.

	We claim that $d_{\bS}(x_i^*, x_j^*) < r$ for all $1 \le i,j \le k$.
	If this is not the case, then there exist indices $i,j$ and $1 \leq a \leq q-1$ such that $d_{\bS}(x_i^*, \omega^a x_j^*) < r$.
	Applying the triangle inequality to $x_i^*, x_j^*$ and $x_0$, we have
	\[d_{\bS}(x_i^*, x_j^*) \leq d_{\bS}(x_i^*, x_0) + d_{\bS}(x_0, x_j^*) < 2r.\]
	On the other hand, applying the triangle inequality to $x_i^*, x_j^*$ and $\omega^a x_j^*$, we have
	\[
	d_{\bS}(x_i^*, x_j^*) \geq d_{\bS}(x_j^*, \omega^a x_j^*) - d_{\bS}(\omega^a x_j^*, x_i^*) > \frac{2\pi}{q} - r,
	\]
	given $1 \leq a \leq q-1$.
	Combining the above two inequalities, since $q \geq 2$, we have
	\[
	2r > \frac{2\pi}{q} - r \implies r > \frac{2\pi}{3q} \geq \frac{2\pi}{q(q+1)},
	\]
	which contradicts the assumption $r \le \frac{2\pi}{q(q+1)}$.
	Thus, $d_{\bS}(x_i^*, x_j^*) < r$ holds for all $1 \le i,j \le k$, completing the proof of the claim.

	We still need to verify the uniqueness of such $\{g_i\}$.
	This is because $x_i^*$ and $gx_i^*$, for $g$ a non-identity element, cannot simultaneously satisfy $d_{\bS}(x_0, x_i^*) < r$ and $d_{\bS}(x_0, gx_i^*) < r$.
	If they were, then we would have
	\[\tfrac{2\pi}{q} \leq d_{\bS}(x_i^*, gx_i^*) \leq d_{\bS}(x_i^*, x_0) + d_{\bS}(x_0, gx_i^*) < 2r \leq \tfrac{2\pi}{q} \cdot \tfrac{2}{q+1},
	\]
	contradicting the assumption $q \geq 2$.
\end{proof}

\subsection{Quotients of spheres}

\subsubsection{}\label{ss:VR-compatible-Sn}

A proper \(G\)-action by isometries on \(\bS^n\) is said to be \defn{\(\VR\)-compatible} if:
\begin{enumerate}
	\item There is \(r_0 > 0\) such that the action is strong \(r\)-diameter for all \(0 < r < r_0\).
	\item For each \(r \in (0, \zeta_n]\), the action commutes with the canonical projection of \(\VR_r(\bS^n)\) and the induced map on orbits
	\[
	\tilde f_r^n \colon \VR_r(\bS^n)_G \to \bS^n_G
	\]
	is a weak equivalence.
\end{enumerate}

For example, based on the results in \cref{ss:VRSn projection}, the trivial action is \(\VR\)-compatible.

\subsubsection{}\label{subsub:VR-compatible-system}

An \defn{equatorial system} is a diagram
\[
\bS^{n_1} \to \bS^{n_2} \to \bS^{n_3} \to \dotsb
\]
of round spheres where each map is an isometric embedding.
For example, the real and complex equatorial systems
\[
\bS^1 \subset \bS^2 \subset \bS^3 \subset \dotsb
\quad\text{and}\quad
\bS^1 \subset \bS^3 \subset \bS^5 \subset \dotsb,
\]
used for the definition of real projective (\cref{sss:cohomology_rpn}) and Lens spaces (\cref{sss:cohomology_lens}), are defined by the canonical inclusions
\[
\R^2 \subset \R^3 \subset \R^4 \subset \dotsb
\quad\text{and}\quad
\bC \subset \bC^2 \subset \bC^3 \subset \dotsb.
\]

\subsubsection{}\label{ss:system VR compatible}

A \(G\)-action on an equatorial system consists of a \(G\)-action on each sphere commuting with the isometric embeddings.
It is said to be free, proper, or by isometries if the \(G\)-action on every sphere is.

The group \(\cyc_2 = \set{1,-1} \subset \R\) acts on the real equatorial system, and, for \(q \in \N\), the multiplicative subgroup \(\cyc_q\) of \(q\)-roots of unity acts on the complex one.
Both of these actions are free, proper, and by isometries.

We say that a proper \(G\)-action by isometries on an equatorial system is \defn{\(\VR\)-compatible} if the \(G\)-action on each sphere is \(\VR\)-compatible and the following diagram commutes for every \(i \in \N\) and $0 < r < \zeta_{n_{i+1}}$:
\begin{equation}\label{eq:VR_quotient}
	\begin{tikzcd}
		\VR_r(\bS^{n_i})_G
		\ar[d]
		\ar[r, "\tilde f_{\,r}^{\,n_i}" above]
		&
		\bS^{n_i}_G
		\ar[d]
		\\
		\VR_r(\bS^{n_{i+1}})_G
		\ar[r, "\tilde f_{\,r}^{\,n_{i+1}}" above]
		&
		\bS^{n_{i+1}}_G.
	\end{tikzcd}
\end{equation}

\lemma The above \(\cyc_2\)-action on the real equatorial system and the \(\cyc_q\)-action on the complex one are \(\VR\)-compatible for every \(q \in \N\).

\begin{proof}
	Let $n \in \N$ and take any $0 < r \leq \zeta_{n+1} (< \zeta_{n})$.
	The commutativity of Diagram~\ref{eq:VR_quotient} for every $i \in \N$ follows directly from the commutativity of the analogous diagram without group actions.

	Next, we verify that the group action on each sphere $\bS^n$ is \(\VR\)-compatible (see \cref{ss:VR-compatible-Sn}).
	By \cref{sss:strong_r_action}, these actions are strong \(r\)-diameter actions for small enough $r>0$.
	We still need to verify that the group action commutes with the $\VR_r(\bS^{n})$-projection $f_r^{n}\colon \VR_r(\bS^{n}) \to \bS^{n}$ and that $f_r^{n}$ induces weak equivalences on the orbit spaces for all $0 < r \leq \zeta_{n}$.

	For the $\cyc_2$-action on the real equatorial system, verifying that the $\cyc_2$-action commutes with $f_r^{n}$ is straightforward.
	By \cref{ss:VRSn projection}, $f_r^{n}$ is a homotopy equivalence.
	Because the \(\cyc_2\)-action is proper and free, the quotient maps to the orbit spaces are covering space projections \cite[Proposition~1.40]{hatcher2000}.
	These projections induce isomorphisms on \(\pi_k\) for \(k \ge 2\) \cite[Proposition~4.1]{hatcher2000}, which, combined with the \(\cyc_2\)-equivariance of $f_r^n$, forces the induced map on the orbit spaces to also be an isomorphism for \(k \ge 2\).
	For \(k=1\), a diagram chase in the short exact sequences of groups of the covering spaces ensures an isomorphism on fundamental groups.
	Together with a natural bijection on path components, this establishes that the induced map on the orbit spaces is a weak equivalence.

	The case of $\cyc_q$-action on the complex equatorial system follows similarly.
\end{proof}

\subsubsection{}\label{ss:fundamental_lemma}

We establish the following lemma, which will be used in \cref{s:critical_radii_rpn} to obtain the homological radii of real projective spaces.
Additionally, we expect that its usefulness will extend to Lens spaces; see \cref{s:critical_radii_lens}.

\lemma Let
\[
\bS^{n_1} \to \bS^{n_2} \to \bS^{n_3} \to \dotsb .
\]
be either the real or complex equatorial system with \(G\) respectively \(\cyc_2\) or \(\cyc_q\) for \(q\) an odd prime.
If $\fillrad(\bS^{n_i}_G)$ is non-decreasing as a function of \(i\) then, for any \(i \leq j\),
\[
\firstdeath{n_i}{\bS^{n_j}_G} \leq \fillrad(\bS^{n_j}_G).
\]

\begin{proof}
	To simplify notation, for any $i \leq j$, let
	\[
	\cX_j = \bS^{n_j}_G, \,\delta_j = \fillrad(\cX_j) \text{ and }\beta_i^j = \firstdeath{n_i}{\cX_j}.
	\]
	We will use an induction argument on $j$.
	When $j = 1$, because $\cX_1$ is connected and $n_1$-dimensional, we apply results in \cref{ss:beta v.s. fillrad} to deduce that $\beta_1^1 = \delta_1$.

	Assume the statement holds for $\cX_{j-1}$, that is, $\beta_i^{j-1} \leq \delta_{j-1}$ for any $i \leq j-1$.
	Since $\delta_{j-1} \leq \delta_j$, we have the following commutative diagram of topological spaces for any $r,\epsilon>0$ small:
	\begin{equation}\label{eq:fundamental_bars_diagram}
		\begin{tikzcd}
			\cX_{j-1}
			\ar[d, hook,"{\iota}" left]
			&
			\VR_r(\cX_{j-1})
			\ar[d, hook,"\iota_r" left]
			\ar[l, "\rho^{j-1}" above]
			\ar[r, hook, "v^{j-1}"]
			&
			\VR_{2\delta_{j-1}+\epsilon}(\cX_{j-1})
			\ar[d, hook]
			\\
			\cX_j
			&
			\VR_r(\cX_j)
			\ar[l, "\rho^j" below]
			\ar[r, hook, "v^j" below]
			&
			\VR_{2\delta_j+\epsilon}(\cX_j).
		\end{tikzcd}
	\end{equation}
	Here, the vertical maps are all induced by the equatorial inclusion of the orbit spaces $\iota \colon \cX_{j-1} \hookrightarrow \cX_j$.
	The horizontal inclusion $v^{j-1}$ (resp. $v^j$) in the right-hand-side square is the inclusion map in the corresponding Vietoris--Rips filtration.
	The horizontal map $\rho^{j-1}$ in the left-hand-side square is the composition of the following maps introduced in \cref{ss:h} and \cref{ss:VRSn projection}, respectively:
	\[\VR_r(\cX_{j-1}) \xrightarrow{h^{-1}_r}\VR_r(\bS^{n_{j-1}})_G \xrightarrow{\tilde{f}_r^{n_{j-1}}} \cX_{j-1}.\]
	Since the \(G\)-action on \(\bS^{n_{j-1}}\) is VR-compatible, it is a strong \(r\)-diameter action for small enough $r>0$, making \(h_r\) an isomorphism by \cref{ss:h} and thus we can consider the inverse map $h^{-1}_r$.
	This VR-compatibility further implies that \(\tilde{f}_r^{n_{j-1}}\) is a weak equivalence, and therefore \(\rho^{j-1}\) is as well.
	The horizontal map \(\rho^j\), defined similarly, is likewise a weak equivalence.

	For any \(i \leq j - 1\), apply the degree \(n_i\) reduced homology functor with appropriate finite field coefficients to the diagram \eqref{eq:fundamental_bars_diagram}.
	Using that the \(G\)-action on the system is \(\VR\)-compatible (cf~\cref{ss:system VR compatible}), we obtain a commutative diagram of vector spaces for $r,\epsilon > 0$ small:
	\begin{equation}\label{eq:diagram of H}
		\begin{tikzcd}[column sep = 4.5em]
			\rH_{n_i}(\cX_{j-1})
			\ar[d, "\cong" left]
			&
			\rH_{n_i}\big(\VR_r(\cX_{j-1})\big)
			\ar[d, "\rH_{n_i}(\iota_r)" left, "\cong" right, myred]
			\ar[l, "\cong" above]
			\ar[r, "\rH_{n_i}(v^{{j-1}})", myred]
			&
			\rH_{n_i}\big(\VR_{2\delta_{j-1}+\epsilon}(\cX_{j-1})\big)
			\ar[d]
			\\
			\rH_{n_i}(\cX_j)
			&
			\rH_{n_i}\big(\VR_r(\cX_j)\big)
			\ar[l, "\cong"]
			\ar[r, "\rH_{n_i}(v^j)" below, myred]
			&
			\rH_{n_i}\big(\VR_{2\delta_j+\epsilon}(\cX_j)\big).
		\end{tikzcd}
	\end{equation}
	The left-most vertical map is an isomorphism in all degrees $\le n_{j-1}$.
	Indeed, with appropriate finite field coefficients, the induced map on cohomology sends algebra generators (see \S\ref{sss:cohomology_rpn} and \S\ref{sss:cohomology_lens}) to algebra generators, hence is an isomorphism in all degrees $\le n_{j-1}$.
	By the universal coefficient theorem, the same holds for homology.

	Let $\sigma_i$ be a representative class for the bar $(0,\, 2\beta_{i}^{j-1})$ in $\VR(\cX_{j-1})$.
	Commutativity of the left-hand-side square of Diagram~(\ref{eq:diagram of H}) implies that $\rH_{n_i}(\iota_r)$ is an isomorphism.
	As $r$ is arbitrarily small, we obtain that $\rH_{n_i}(\iota_r)(\sigma_i)$ creates a bar born at $0$.
	Moreover, this bar dies before or at $2\delta_j$ as explained below.
	It follows from the induction hypothesis $\beta_i^{j-1} \leq \delta_{j-1}$ that $\sigma_i$ dies before $2\delta_{j-1} + \epsilon$, implying that $\rH_{n_i}(v^{{j-1}})(\sigma_i) = 0$.
	Using the right-hand-side square's commutativity, we have $(\rH_{n_i}(v^j) \circ \rH_{n_i}(\iota_r))(\sigma_i)=0$, which means $\rH_{n_i}(\iota_r)(\sigma_i)$ dies before $2\delta_j+\epsilon$.
	Since this holds for any \(\epsilon > 0\), we conclude that \(\beta_i^j \leq \delta_j\).

	For the case when $i = j$, apply results in \cref{ss:beta v.s. fillrad} again to get that $\beta_j^j = \delta_j$.
	This completes the proof.
\end{proof}


\subsection{Real projective spaces}\label{s:critical_radii_rpn}

Let us now focus on the real equatorial system with its canonical \(\rC_2\)-action.
We will consider every sphere \(\bS^n(2)\) to have radius 2 so the quotient \(\rp^n\) has the same diameter as \(\bS^n\), the round sphere of radius 1.
The ground field here is \(\Ftwo\).

\lemma
For \(k,n,m \in \N\) and \(\Sq^k \in \cO(m-k, m)\) we have:
\begin{enumerate}
	\item \(\crit(\rp^n) = \frac{\pi}{3}\).
	\item \(\firstdeath{m}{\rp^n} =
	\begin{cases}
		\frac{\pi}{3} & m \leq n, \\
		\hfil 0 & \text{otherwise}.
	\end{cases}\)
	\item \(\firstdeath{\Sq^k}{\rp^n} =
	\begin{cases}
		\tfrac{\pi}{3} & m \leq n \text{ and } \binom{m-k}{k} \text{ is odd},\\
		\hfil 0 & \text{otherwise}. \\
	\end{cases}\)
\end{enumerate}

\begin{proof}
	(1) By \cite{katz1983filling}, the filling radius of $\rp^n$ is $\frac{\pi}{3}$.
	Since the homotopical critical radius is always bounded above by the filling radius, we have $\crit(\rp^n) \leq \frac{\pi}{3}$.
	The reverse inequality, $\crit(\rp^n) \geq \frac{\pi}{3}$, follows from the remark following \cite[Theorem~4.5]{adams2022metric}.

	\smallskip (2) We only need to prove for the case of $m\leq n$, since the other case follows from the definition since \(\rH_m(\rp^n) = 0\) if \(m > n\) or \(m = 0\).
	Consider the \(\rC_2\)-action on the real equatorial system \(\bS^{1} \to \bS^{2} \to \dotsb \to \bS^{n}\).
	By \cref{ss:system VR compatible}, this system is \(\VR\)-compatible.
	By \cite{katz1983filling}, the filling radius of the quotient space $\bS^{i}_{\rC_2}$ is $\frac{\pi}{6}$ for any $i \geq 1$, and thus $\fillrad(\bS^{i}_{\rC_2})$ is non-decreasing as a function of \(i\).
	Therefore, we can apply \cref{ss:fundamental_lemma} to deduce that for any \(m \leq n\),
	\[
	\firstdeath{m}{\bS^{n}_{\rC_2}} \leq \fillrad(\bS^{n}_{\rC_2}) = \tfrac{\pi}{6}.
	\]
	On the other hand, $\firstdeath{m}{\bS^{n}_{\rC_2}}\geq \crit(\bS^{n}_{\rC_2}) = \tfrac{\pi}{6}$.
	Thus, the equality holds.
	As $\rp^n$ is the quotient of the $n$-sphere with radius $2$,
	\[\firstdeath{m}{\rp^n} = \firstdeath{m}{\bS^{n}(2)_{\rC_2}} = 2\cdot \firstdeath{m}{\bS^{n}_{\rC_2}} = \frac{\pi}{3}.\]

	\smallskip (3) If $m > n$ or $\binom{m-k}{k}$ is even, then \(\img((\Sq^k)_{\rp^n}) = 0\).
	Assume therefore that $m \leq n$ and $\binom{m-k}{k}$ is odd.
	Because $2\crit(\rp^n)=\frac{2\pi}{3}$, $\VR_r(\rp^n)$ retains the homotopy type of $\rp^n$ for $r \in (0,\tfrac{2\pi}{3})$.
	Recall from \cref{sss:cohomology_rpn} the cohomology algebra \(\rH^\ast(\rp^n; \Ftwo) \cong \frac{\Ftwo[\sigma]}{(\sigma^{n+1})}.\)
	Thus, $\Sq^k(\sigma^{m-k}) = \sigma^{m}$ generates a bar in the $\img_{\Sq^k}$-barcode that is born at $0$ and stays alive until the non-trivial class $\sigma^{m}$ dies at $\tfrac{2\pi}{3}$.
	Thus, $\firstdeath{\Sq^k}{\rp^n} \leq \tfrac{\pi}{3}$.
	On the other hand, $\firstdeath{\Sq^k}{\rp^n} \geq \crit(\rp^n) = \tfrac{\pi}{3}$.
\end{proof}

Using the above values, the estimates resulting from the analysis of \cref{ss:barcode_general} are illustrated in \cref{fig:sq barcodes}.

\begin{figure}
	\centering
	\begin{tikzpicture}[scale=0.52]
	\begin{axis} [
		title = {\LARGE $\Hbarc{m}{\rp^n},\, m\leq n$},
		ticklabel style = {font=\Large},
		axis y line=middle,
		axis x line=middle,
		ytick={0.5,0.67,0.95},
		yticklabels={$\frac{\pi}{2}$,$\frac{2\pi}{3}$,$\pi$},
		xtick={0.5,0.67,0.95},
		xticklabels={$\frac{\pi}{2}$,$\frac{2\pi}{3}$,$\pi$},
		xmin=-0.015, xmax=1.1,
		ymin=0, ymax=1.1,]
		\addplot [mark=none] coordinates {(0,0) (1,1)};
		\addplot [thick,color=black!20!white,fill=black!30!white,
		fill opacity=0.4]coordinates {
			(0.67,0.95)
			(0.67,0.67)
			(0.95,0.95)
			(0.67,0.95)};
		\addplot [black!40!white,mark=none,dashed, thin] coordinates {(0,0.67) (0.67,0.67)};
		\addplot [black!40!white,mark=none,dashed, thin] coordinates {(0.67,0) (0.67,0.67)};
		\addplot[barccolor,mark=*] (0, 0.67) circle (2pt) node[above right,barccolor]{};
	\end{axis}
\end{tikzpicture}
\begin{tikzpicture}[scale=0.52]
	\begin{axis} [
		title={\LARGE $\Hbarc{m}{\rp^n},\, m>n$},
		ticklabel style = {font=\Large},
		axis y line=middle,
		axis x line=middle,
		ytick={0.5,0.67,0.95},
		yticklabels={$\frac{\pi}{2}$,$\frac{2\pi}{3}$,$\pi$},
		xtick={0.5,0.67,0.95},
		xticklabels={$\frac{\pi}{2}$,$\frac{2\pi}{3}$,$\pi$},
		xmin=-0.015, xmax=1.1,
		ymin=0, ymax=1.1,]
		\addplot [mark=none] coordinates {(0,0) (1,1)};
		\addplot [thick,color=black!20!white,fill=black!30!white,
		fill opacity=0.4]coordinates {
			(0.67,0.95)
			(0.67,0.67)
			(0.95,0.95)
			(0.67,0.95)};
		\addplot [black!40!white,mark=none,dashed, thin] coordinates {(0,0.67) (0.67,0.67)};
		\addplot [black!40!white,mark=none,dashed, thin] coordinates {(0.67,0) (0.67,0.67)};
	\end{axis}
\end{tikzpicture}

\begin{tikzpicture}[scale=0.52]
	\begin{axis} [
		title = {\LARGE $\sqbarcl{k}{}{\rp^n},\, m \leq n$ and $\binom{m-k}{k}$ odd},
		ticklabel style = {font=\Large},
		axis y line=middle,
		axis x line=middle,
		ytick={0.5,0.67,0.95},
		yticklabels={$\frac{\pi}{2}$,$\frac{2\pi}{3}$,$\pi$},
		xtick={0.5,0.67,0.95},
		xticklabels={$\frac{\pi}{2}$,$\frac{2\pi}{3}$,$\pi$},
		xmin=-0.015, xmax=1.1,
		ymin=0, ymax=1.1,]
		\addplot [mark=none] coordinates {(0,0) (1,1)};
		\addplot [thick,color=black!20!white,fill=black!30!white,
		fill opacity=0.4]coordinates {
			(0.67,0.95)
			(0.67,0.67)
			(0.95,0.95)
			(0.67,0.95)};
		\addplot [black!40!white,mark=none,dashed, thin] coordinates {(0,0.67) (0.67,0.67)};
		\addplot [black!40!white,mark=none,dashed, thin] coordinates {(0.67,0) (0.67,0.67)};
		\addplot[barccolor,mark=*] (0, 0.67) circle (2pt) node[above right,barccolor]{};
	\end{axis}
\end{tikzpicture}
\begin{tikzpicture}[scale=0.52]
	\begin{axis} [
		title={\LARGE $\sqbarcl{k}{}{\rp^n}$, otherwise},
		ticklabel style = {font=\Large},
		axis y line=middle,
		axis x line=middle,
		ytick={0.5,0.67,0.95},
		yticklabels={$\frac{\pi}{2}$,$\frac{2\pi}{3}$,$\pi$},
		xtick={0.5,0.67,0.95},
		xticklabels={$\frac{\pi}{2}$,$\frac{2\pi}{3}$,$\pi$},
		xmin=-0.015, xmax=1.1,
		ymin=0, ymax=1.1,]
		\addplot [mark=none] coordinates {(0,0) (1,1)};
		\addplot [thick,color=black!20!white,fill=black!30!white,
		fill opacity=0.4]coordinates {
			(0.67,0.95)
			(0.67,0.67)
			(0.95,0.95)
			(0.67,0.95)};
		\addplot [black!40!white,mark=none,dashed, thin] coordinates {(0,0.67) (0.67,0.67)};
		\addplot [black!40!white,mark=none,dashed, thin] coordinates {(0.67,0) (0.67,0.67)};
	\end{axis}
\end{tikzpicture}
	\caption{Let $m \in \N$ and $\Sq^k \in \cO(m-k, m)$.
		\emph{Top row:} the persistent reduced homology barcode of $\rp^n$.
		\emph{Bottom row:} the $\img_{\Sq^k}$-barcode of $\rp^n$.
		In each figure, the gray region represents where additional bars could exist within the corresponding barcode.}
	\label{fig:sq barcodes}
\end{figure}

\subsection{Lens spaces}\label{s:critical_radii_lens}

Much less information is known about the critical radii of Lens spaces (\cref{sss:cohomology_lens}).
In this section, we establish a relationship between the filling radius and the homological radii of Lens spaces, assuming a monotonicity condition on the former.
For $n \geq 1$, the quotient of the unit sphere \(\bS^{2n+1}\) by \(\rC_q\), with $q$ an odd prime, has diameter \(\frac{\pi}{2}\).
We equip every \(\rL_q^n\) with the quotient metric induced by \(\bS^{2n+1}(2)\); for $n \geq 1$, this gives \(\diam(\rL_q^n) = \pi\).
The ground field here is $\bF_q$.

\lemma
For fixed $n \geq 1$, if $\fillrad(\rL^0_q) \leq \dotsb \leq \fillrad(\rL^n_q)$, then for any odd degree $m$ with $1 \leq m \leq 2n+1$,
\[
\firstdeath{m}{\rL^n_q} \leq \fillrad(\rL^n_q).
\]

\begin{proof}
	We apply an argument similar to that in the second part of the proof of \cref{s:critical_radii_rpn}.
	By \cref{ss:system VR compatible}, the $\rC_q$-action on the system $\bS^1 \subset \bS^3 \subset \dotsb \subset \bS^{2n+1}$ of unit spheres is \(\VR\)-compatible.
	This, combined with the assumption that $\fillrad\big(\bS^{2\tilde{n}+1}_{\rC_q}\big) = \tfrac{1}{2} \fillrad(\rL^{\tilde{n}}_q)$ is non-decreasing in $\tilde{n}$ for all $0\leq \tilde{n} \leq n$, meets the conditions of \cref{ss:fundamental_lemma}, implying that for any odd degree $m$ with $1 \leq m \leq 2n+1$,
	\[
	\mathrm{Rad}_m\big(\bS^{2n+1}_{\rC_q}\big) \leq \fillrad\big(\bS^{2n+1}_{\rC_q}\big)
	\]
	Consequently, for any such odd degree $m$, we obtain
	\[
	\firstdeath{m}{\rL^n_q} \leq \fillrad(\rL^n_q). \qedhere
	\]
\end{proof}


\section{Gromov--Hausdorff estimates}\label{s:gh_estimates}

In this section, we construct examples of Riemannian pseudomanifolds where Gromov--Hausdorff estimates from persistent cohomology operations yield tighter bounds than those from persistent homology.
This demonstrates the enhanced discriminating power of these persistent invariants in Riemannian geometry.

\subsection{General estimates}\label{ss:genberal_distance_comparison}

\subsubsection{}

We will prove the relevant statements in greater generality, allowing for similar estimates to be established for Lens spaces once their critical radii are determined.

\medskip\noindent Let $\cM$ be a closed connected Riemannian \(n\)-manifold of diameter $\pi$.

\medskip\noindent\textit{\darkblue Desiderata}.\
\begin{enumerate}
	\item $\crit(\cM) \geq \tfrac{\zeta_n}{2} $;
	\item $\firstdeath{m}{\cM} < \tfrac{\zeta_n}{4}+\tfrac{\zeta_m}{2}$ for any $m \in \N$.
	\item There exists $\theta \in \cO(\ell, m)$ with \(\ell \neq m\) such that
	\[
	\img_\theta(\cM) \neq 0
	\quad\text{and}\quad
	\dim \opH_m(\cM) = 1.
	\]
\end{enumerate}

\medskip We verify that the round real projective spaces of diameter \(\pi\) satisfy these Desiderata.
Recall that \(\zeta_n = \arccos(\tfrac{-1}{1+n})\) and, from \cref{s:critical_radii_rpn}, that \(\crit(\rp^n) = \tfrac{\pi}{3}\).
For \(n \geq 1\) we have \(\crit(\rp^n) \geq \tfrac{\zeta_n}{2}\) since
\[
\tfrac{1}{2} \geq \tfrac{1}{1+n} \ \Leftrightarrow \
\tfrac{-1}{2} \leq \tfrac{-1}{1+n} \ \Leftrightarrow \
\arccos \tfrac{-1}{2} \geq \arccos\tfrac{-1}{1+n} \ \Leftrightarrow \
\tfrac{2\pi}{3} \geq \zeta_n.
\]

Let us now focus on Desideratum~(2).
For any integers $k \geq 1$ we have $\frac{-1}{1+k} < 0$.
Because \(\arccos\) is decreasing, it follows that for any $k \geq 1$,
\[
\zeta_k = \arccos\left(\frac{-1}{1+k}\right) > \arccos(0) = \frac{\pi}{2}.
\]
Using this strict lower bound for both $\zeta_n$ and $\zeta_m$, we obtain:
\[
\frac{\zeta_n}{4} + \frac{\zeta_m}{2} > \frac{1}{4}\left(\frac{\pi}{2}\right) + \frac{1}{2}\left(\frac{\pi}{2}\right) = \frac{\pi}{8} + \frac{\pi}{4} = \frac{3\pi}{8} > \frac{\pi}{3}.
\]
By \cref{s:critical_radii_rpn} we know that either \(\firstdeath{m}{\rp^n} = \tfrac{\pi}{3}\) or \(\firstdeath{m}{\rp^n} = 0\) from which the above inequality concludes the verification of Desideratum~(2).

As shown in \cref{s:critical_radii_rpn}, for every $n>1$ Desideratum~(3) is satisfied by
\[
\Sq^1 \colon \rH^1(\rp^n;\Ftwo) \to \rH^2(\rp^n;\Ftwo),
\]
since \(\Sq^1(\sigma)=\sigma^2\neq 0\) and \(\dim \opH_2(\rp^n;\Ftwo)=1\).

\medskip For fixed $n \geq 1$ and the round Lens space $\rL^n_p$ of diameter \(\pi\), Desideratum~(3) is achieved by selecting $\theta$ as the Bockstein \(\beta \colon \rH^1(-;\Fp) \to \rH^2(-;\Fp)\) (see \cref{sss:cohomology_lens}); its target homology group is one-dimensional.
Desideratum~(1) remains unproven.
For odd degrees, Desideratum~(2) follows if
\[
\fillrad(\rL^0_p) \leq \dotsb \leq \fillrad(\rL^n_p) < \tfrac{3\zeta_{2n+1}}{4}.
\]
Indeed, under this condition, \cref{s:critical_radii_lens} implies that for every odd degree $m \leq 2n+1$,
\[
\firstdeath{m}{\rL^n_p} \leq \fillrad(\rL^n_p) < \tfrac{3\zeta_{2n+1}}{4} \leq \tfrac{\zeta_{2n+1}}{4} + \tfrac{\zeta_{m}}{2}.
\]
For even degrees, the corresponding inequalities in Desideratum~(2) remain additional conditions to be verified.

\subsubsection{}\label{sss:choice of S_M}

Let \(\cM\) be a closed connected Riemannian \(n\)-manifold of diameter \(\pi\) satisfying our Desiderata.
Let $\bS_{\cM}$ be the wedge sum of spheres of diameter \(\pi\) with the same reduced homology groups as $\cM$.

\theorem
There is $\ell, m \in \N$ and $\theta \in \cO(\ell, m)$ such that
\begin{equation}\label{eq:comparison}
	\di \Big(\rH^{\VR}_{m}(\bS_{\cM}),\, \rH^{\VR}_{m}(\cM)\Big) < \di \big(\img^{\VR}_{\theta}(\bS_{\cM}),\, \img^{\VR}_{\theta}(\cM)\big).
\end{equation}

The proof of this theorem involves the following three steps, which we will respectively prove in \cref{sss:db_upper_bound}, \cref{sss:db_theta_lower_bound}, and \cref{sss:comparison_lemma_zeta_n}.
\begin{enumerate}
	\item [(a)] LHS of \cref{eq:comparison} is bounded above by $\max\big\{|\zeta_m - 2\firstdeath{m}{\cM} |, \tfrac{\pi - \zeta_n}{2}\big\}$.
	\item [(b)] RHS of \cref{eq:comparison} is bounded below by $\min\{\firstdeath{\theta}{\cM}, \zeta_n\}$.
	\item [(c)] The value in (a) is strictly less than the value in (b).
\end{enumerate}

\subsubsection{}\label{sss:db_upper_bound}

\lemma
If \(\dim \opH_m(\cM) = 1\), then
\begin{equation}\label{eq:db_usual_upper_bound}
	\di\Big(\rH^{\VR}_{m}(\bS_{\cM}),\, \rH^{\VR}_{m}(\cM)\Big)
	\leq \max\big\{|\zeta_m - 2\firstdeath{m}{\cM} |,\, \tfrac{\pi - \zeta_n}{2}\big\}.
\end{equation}
If \(\opH_m(\cM)\) is zero, a similar result holds, with the upper bound being $\tfrac{\pi - \zeta_n}{2}$.

\begin{proof}
	Because both metric spaces are totally bounded, their persistent homology is q-tame.
	By \cref{ss:algebraic_stability}, the interleaving distance between their persistent homology equals the bottleneck distance between the corresponding barcodes.

	To obtain an upper bound on the bottleneck distance, it suffices to construct a matching and calculate its cost.
	Let \(\beta_m = 2\firstdeath{m}{\cM}\).
	First, suppose \(\dim \opH_m(\cM) = 1\).
	Consider the matching \(P\) between \(\Hbarc[\field]{m}{\bS_{\cM}}\) and \(\Hbarc[\field]{m}{\cM}\) such that the bar \((0, \zeta_m)\) from \(\Hbarc[\field]{m}{\bS_{\cM}}\) is matched to \((0, \beta_m)\) from \(\Hbarc[\field]{m}{\cM}\), with all other bars remaining unmatched. Using the estimates from \cref{fig:barcodes_general} and \cref{fig:barcodes_vs}, we find that the cost of \(P\) satisfies
	\begin{align*}
		\cost(P)
		\leq & \max\Big\{|\zeta_m - \beta_m|,\, \frac{\pi - \zeta_n}{2},\, \frac{\pi - 2\crit(\cM)}{2}\Big\} \\
		= & \max\Big\{|\zeta_m - \beta_m|,\, \frac{\pi - \zeta_n}{2}\Big\} \quad \text{(since \(\zeta_n \leq 2\crit(\cM)\))}.
	\end{align*}
	Thus, \cref{eq:db_usual_upper_bound} follows.

	Now, suppose \(\opH_m(\cM)\) is zero. In this case, the matching where all bars are unmatched has a cost no larger than \(\frac{\pi - \zeta_n}{2}\).
\end{proof}

\subsubsection{}\label{sss:db_theta_lower_bound}

\lemma
Let $\theta \in \cO(\ell, m)$ with \(\ell \neq m\) be such that $\img \theta_\cM$ is non-zero.
Then,
\begin{equation}\label{eq:db_theta_lower_bound}
	\di\big(\img^{\VR}_{\theta}(\bS_{\cM}), \img^{\VR}_{\theta}(\cM)\big)
	\geq \min\{\firstdeath{\theta}{\cM}, \zeta_n\}.
\end{equation}

\begin{proof}
	Since both metric spaces are totally bounded, their persistent cohomology is q-tame, making the persistent $\img_\theta$-modules q-tame as well (see \cref{ss:theta-modules-q-tame}).
	Thus, the left-hand-side of \cref{eq:db_theta_lower_bound} can be replaced with the bottleneck distance between the corresponding $\img_\theta$-barcodes.

	To establish a lower bound on this bottleneck distance, we need to compute the minimum cost of matching specific bars.

	Let \(\gamma_\theta = 2\firstdeath{\theta}{\cM}\).
	From \cref{fig:barcodes_general}, note that \(\thetabarc{\cM}\) contains a bar \((0, \gamma_\theta)\).
	We now calculate the minimum cost of matching this bar.

	Consider an arbitrary matching \(Q\) between \(\thetabarc{\bS_{\cM}}\) and \(\thetabarc{\cM}\).
	There are two cases to analyze.
	If $(0,\gamma_\theta )$ is unmatched in $Q$, then $\cost(Q)$ is at least $ \tfrac{\gamma_\theta }{2}$.
	If $(0,\gamma_\theta )$ is matched to some bar $(a,b) \in \thetabarc{\bS_{\cM}}$, then
	$\cost(Q) \geq \|(0,\gamma_\theta ) - (a,b)\|_\infty \geq a \geq \zeta_n$, because $(a,b) \subset (\zeta_n, \pi)$.
	Thus, any matching $Q$ must satisfy $\cost(Q) \geq \min\{\tfrac{\gamma_\theta }{2}, \zeta_n\}$, which proves \cref{eq:db_theta_lower_bound}.
\end{proof}

\subsubsection{}\label{sss:comparison_lemma}

The following results will be used in \cref{sss:comparison_lemma_zeta_n} and \cref{ss:distance_estimate_rpn}.

As before, let $\zeta_n = \arccos\!\left(-\tfrac{1}{n+1}\right)$.
Since $\arccos$ is strictly decreasing on $[-1,1]$ and $-\frac{1}{n+1}$ is increasing in $n$, the sequence $(\zeta_n)_{n\ge 1}$ is strictly decreasing.

\lemma
Let $m \leq n$.
\begin{enumerate}
	\item $\zeta_m - \zeta_n < \tfrac{\pi - \zeta_n}{2} < \tfrac{\zeta_n}{2}$.
	\item For any $r \in \big[\zeta_n,\, \tfrac{\zeta_n}{2}+\zeta_m\big)$, $\max\big\{|\zeta_m - r|,\, \tfrac{\pi - \zeta_n}{2}\big\} < \tfrac{\zeta_n}{2}$.
\end{enumerate}

\begin{proof}
	Since $-\frac{1}{n+1} < 0$, we have $\zeta_n > \arccos(0) = \tfrac{\pi}{2}$.
	Since $(\zeta_n)$ is decreasing and $\zeta_1 = \arccos\!\left(-\tfrac12\right) = \tfrac{2\pi}{3}$, we also have $\zeta_n \le \tfrac{2\pi}{3}$.

	\smallskip
	(1) Because $(\zeta_n)$ is decreasing and $m\le n$, we have $\zeta_m \le \zeta_1$, hence
	$\zeta_m - \zeta_n \le \zeta_1 - \zeta_n < \tfrac{2\pi}{3} - \tfrac{\pi}{2} = \tfrac{\pi}{6}$.
	Since $\zeta_n \le \tfrac{2\pi}{3}$, we have $\tfrac{\pi - \zeta_n}{2} \ge \tfrac{\pi}{6}$, so $\zeta_m - \zeta_n < \tfrac{\pi - \zeta_n}{2}$.
	Moreover $\tfrac{\pi - \zeta_n}{2} < \tfrac{\zeta_n}{2}$ is equivalent to $\zeta_n > \tfrac{\pi}{2}$, which holds.

	\smallskip
	(2) Since $\tfrac{\pi - \zeta_n}{2} < \tfrac{\zeta_n}{2}$ by Part (1), it remains to show $|\zeta_m - r| < \tfrac{\zeta_n}{2}$.
	If $r \in [\zeta_n,\, \zeta_m)$, then $|\zeta_m - r| = \zeta_m - r \le \zeta_m - \zeta_n < \tfrac{\zeta_n}{2}$ by Part (1).
	If $r \in [\zeta_m,\, \zeta_m + \tfrac{\zeta_n}{2})$, then $|\zeta_m - r| = r - \zeta_m < \tfrac{\zeta_n}{2}$ by the upper bound on $r$.
	Thus in all cases $|\zeta_m - r| < \tfrac{\zeta_n}{2}$, and therefore $\max\big\{|\zeta_m - r|,\, \tfrac{\pi - \zeta_n}{2}\big\} < \tfrac{\zeta_n}{2}$.
\end{proof}

\subsubsection{}\label{sss:comparison_lemma_zeta_n}

We will finish the proof of \cref{sss:choice of S_M} by completing the third step outlined therein.

\lemma
For any $m \in \N$ and $\theta \in \cO(\ell, m)$ such that $\img \theta_\cM$ is non-zero,
\[\max\big\{|\zeta_m - 2\firstdeath{m}{\cM} |, \tfrac{\pi - \zeta_n}{2}\big\} < \min\{\firstdeath{\theta}{\cM}, \zeta_n\}.\]
\begin{proof}
	Let $\beta_m = 2\firstdeath{m}{\cM}$ and $\gamma_\theta = 2\firstdeath{\theta}{\cM}$.
	Because \(\zeta_n \leq 2\crit(\cM) \leq \beta_m < \tfrac{\zeta_n}{2} + \zeta_m\), $\beta_m$ satisfies the condition in Part (2) of the above lemma.
	Thus, we have the leftmost inequality below:
	\[
	\max\big\{|\zeta_m - \beta_m|, \tfrac{\pi - \zeta_n}{2}\big\} < \tfrac{\zeta_n}{2}
	\leq
	\min\big\{\tfrac{\gamma_\theta}{2}, \zeta_n\big\}.
	\]
	To obtain the rightmost inequality, note that \(\zeta_n \leq 2\crit(\cM) \leq \gamma_\theta\).
\end{proof}

\subsection{Real projective spaces}\label{ss:distance_estimate_rpn}

As stated earlier, the round \(\rp^n\) of diameter \(\pi\) satisfies our Desiderata, so the discussions in the previous subsection apply.
Here we obtain a more precise estimate of the interleaving distances involved.
Let \(\bS_{\rp^n} = \bS^1 \vee \dots \vee \bS^n\) as before.
Then, we have the following result.

\theorem
(1) For any \(n,m \in \N\),
\[
\di\Big(\rH^\VR_m(\bS_{\rp^n}),\, \rH^\VR_m(\rp^n)\Big) \leq \frac{\pi - \zeta_n}{2} < \frac{\pi}{4}.
\]

\noindent (2) For any \(n,m,k \in \N\) such that \(0 < k \leq m \leq n\) and \(\binom{m-k}{k}\) is odd:
\[
\di\Big(\img_{\Sq^k}^\VR(\bS_{\rp^n}),\, \img_{\Sq^k}^\VR(\rp^n)\Big) \geq \frac{\pi}{3},
\]
where $\Sq^k$ lands in degree \(m\).

\begin{proof}
	(1)
	As in the proof of \cref{sss:db_upper_bound}, we first invoke \cref{ss:algebraic_stability} to replace the interleaving distance with the bottleneck distance between the corresponding barcodes, since both metric spaces are totally bounded and their persistent homology is therefore q-tame.

	When $1 \leq m \leq n$, we have shown in \cref{s:critical_radii_rpn} that
	$2\firstdeath{m}{\rp^n} = \tfrac{2\pi}{3}$.
	Thus, the upper bound we computed in \cref{eq:db_usual_upper_bound} on the interleaving distance between persistent homology evaluates as
	\[\max\big\{|\zeta_m - 2\firstdeath{m}{\cM} |, \tfrac{\pi - \zeta_n}{2}\big\} =
	\max \big\{ |\zeta_m - \tfrac{2\pi}{3} |, \tfrac{\pi - \zeta_n}{2} \big\}
	= \tfrac{\pi - \zeta_n}{2}.
	\]
	Here, we applied \cref{sss:comparison_lemma} to deduce that $|\zeta_m - \tfrac{2\pi}{3}| = \zeta_1 - \zeta_m < \tfrac{\pi - \zeta_1}{2} \leq \tfrac{\pi - \zeta_n}{2}$.

	When $m > n$, recall from \cref{fig:barcodes_vs} that bars in the $m$-th persistent homology of $\bS_{\rp^n}$ (if any) are contained in $(\zeta_n, \pi)$, and from \cref{fig:sq barcodes} that bars in the $m$-th persistent homology of $\rp^n$ (if any) are contained in $(\tfrac{2\pi}{3}, \pi)$.
	Thus, the interleaving distance between these persistent homology groups is at most half the maximum bar length (with all bars unmatched), which is no larger than $\max\{\tfrac{\pi - \zeta_n}{2}, \tfrac{\pi}{6}\} \leq \tfrac{\pi - \zeta_n}{2}$ since $\zeta_n \leq \tfrac{2\pi}{3}$.

	(2) With $n,m,k,\Sq^k$ satisfying the required conditions, we have shown in \cref{s:critical_radii_rpn} that
	$\firstdeath{\Sq^k}{\rp^n} = \tfrac{\pi}{3}$.
	Applying \cref{eq:db_theta_lower_bound}, we have
	\[\di\Big(\img_{\Sq^k}^\VR(\bS_{\rp^n}),\, \img_{\Sq^k}^\VR(\rp^n)\Big)
	\geq \min\Big\{\tfrac{\pi}{3},\, \zeta_n\Big\}
	= \tfrac{\pi}{3}.\qedhere
	\]
\end{proof}

By combining Item~(2) of the above theorem with \cref{thm:stability intro}, we obtain a lower bound of \(\frac{1}{2}\cdot\frac{\pi}{3}\) for the Gromov--Hausdorff distance between \(\rp^n\) and \(\bS_{\rp^n}\).
For comparison, applying Item~(1) together with the stability theorem for persistent homology (see page~\pageref{thm:stability ph}) yields a lower bound of at most \(\frac{1}{2}\cdot\frac{\pi}{4}\).
This result demonstrates the stronger distinguishing power of persistent cohomology operations compared to persistent homology.
We do not claim that the bound \(\frac{\pi}{6}\) surpasses all known lower bounds.
For example, the diameter distance \(\frac{1}{2} \lvert \operatorname{diam}(\rp^n) - \operatorname{diam}(\bS_{\rp^n}) \rvert = \frac{1}{2} \lvert \pi - 2\pi \rvert = \frac{\pi}{2}\) is greater.

	\sloppy
	\printbibliography

@article{medina2022per_st,
	author = {{Lupo}, Umberto and {Medina-Mardones}, Anibal M. and {Tauzin}, Guillaume},
	title = "{Persistence Steenrod modules}",
	fjournal = {Journal of Applied and Computational Topology},
	journal = {J. Appl. Comput. Topol.},
	volume={6},
	number={4},
	pages={475--502},
	year = {2022},
	doi = {10.1007/s41468-022-00093-7},
	URL = {https://doi.org/10.1007/s41468-022-00093-7},
	publisher = {Springer}
}

@article{medina2023fast_sq,
	AUTHOR = {{Medina-Mardones}, Anibal M.},
	TITLE = {New formulas for cup-{$i$} products and fast computation of {S}teenrod squares},
	JOURNAL = {Comput. Geom.},
	FJOURNAL = {Computational Geometry. Theory and Applications},
	VOLUME = {109},
	YEAR = {2023},
	ISSN = {0925-7721},
	MRNUMBER = {4473678},
	DOI = {10.1016/j.comgeo.2022.101921},
	URL = {https://doi.org/10.1016/j.comgeo.2022.101921},
}

@article{medina2023fuct_top,
	author = {{Bauer}, Ulrich and {Medina-Mardones}, Anibal M. and {Schmahl}, Maximilian},
	title = "{Persistent homology for functionals}",
	journal = {Commun. Contemp. Math.},
	fjournal = {Communications in Contemporary Mathematics},
	year = {2023},
	doi = {10.1142/S0219199723500554},
	URL = {https://doi.org/10.1142/S0219199723500554},
}

@book {steenrod1962cohomology,
	ISBN = {9780691079240},
	author = {Steenrod, {N. E.} and Epstein, {D. B. A.}},
	publisher = {Princeton University Press},
	title = {Cohomology Operations: Lectures by N. E. Steenrod.},
	year = {1962},
	URL = {http://www.jstor.org/stable/j.ctt1b7x52h},
}

@inproceedings {may1970general,
	AUTHOR = {May, J. Peter},
	TITLE = {A general algebraic approach to {S}teenrod operations},
	BOOKTITLE = {The {S}teenrod {A}lgebra and its {A}pplications},
	SERIES = {Lecture Notes in Mathematics, Vol. 168},
	PAGES = {153--231},
	PUBLISHER = {Springer, Berlin},
	YEAR = {1970},
	MRCLASS = {55.34 (18.00)},
	MRNUMBER = {0281196},
	MRREVIEWER = {W. D. Barcus},
	url = {https://link.springer.com/chapter/10.1007/BFb0058524}
}

@article{lee2017quantifying,
	title={Quantifying similarity of pore-geometry in nanoporous materials},
	author={Lee, Yongjin and Barthel, Senja D and D{\l}otko, Pawe{\l} and Moosavi, S Mohamad and Hess, Kathryn and Smit, Berend},
	journal={Nature communications},
	volume={8},
	number={1},
	pages={1--8},
	year={2017},
	publisher={Nature Publishing Group},
	url={https://doi.org/10.1038/ncomms15396}
}

@article{adams2022gromov,
	title={{G}romov--{H}ausdorff distances, {B}orsuk--{U}lam theorems, and {V}ietoris--{R}ips complexes},
	author={Adams, Henry and Bush, Johnathan and Clause, Nate and Frick, Florian and G{\'o}mez, Mario and Harrison, Michael and Jeffs, R Amzi and Lagoda, Evgeniya and Lim, Sunhyuk and M{\'e}moli, Facundo and others},
	journal={arXiv preprint arXiv:2301.00246},
	url={https://arxiv.org/abs/2301.00246},
	year={2022}
}

@article{ji2021gromov,
	title = {Gromov–Hausdorff distance between interval and circle},
	journal = {Topology and its Applications},
	volume = {307},
	pages = {107938},
	year = {2022},
	issn = {0166-8641},
	doi = {https://doi.org/10.1016/j.topol.2021.107938},
	url = {https://www.sciencedirect.com/science/article/pii/S0166864121003552},
	author = {Yibo Ji and Alexey A. Tuzhilin}
}

@article{lim2021gromov,
	AUTHOR = {Lim, Sunhyuk and M\'{e}moli, Facundo and Smith, Zane},
	TITLE = {The {G}romov-{H}ausdorff distance between spheres},
	JOURNAL = {Geom. Topol.},
	FJOURNAL = {Geometry \& Topology},
	VOLUME = {27},
	YEAR = {2023},
	NUMBER = {9},
	PAGES = {3733--3800},
	ISSN = {1465-3060},
	MRCLASS = {53C23},
	MRNUMBER = {4674839},
	MRREVIEWER = {Wei Zhao},
	DOI = {10.2140/gt.2023.27.3733},
	URL = {https://doi.org/10.2140/gt.2023.27.3733},
}

@article{harrison2023quantitative,
	title={Quantitative upper bounds on the {G}romov-{H}ausdorff distance between spheres},
	author={Harrison, Michael and Jeffs, R Amzi},
	journal={arXiv preprint arXiv:2309.11237},
	url={https://arxiv.org/abs/2309.11237},
	year={2023}
}

@article {saul2024gromov,
	AUTHOR = {Rodr\'{i}guez Mart\'{i}n, Sa\'{u}l},
	TITLE = {Gromov-{H}ausdorff distances from simply connected geodesic
	spaces to the circle},
	JOURNAL = {Proc. Amer. Math. Soc. Ser. B},
	FJOURNAL = {Proceedings of the American Mathematical Society. Series B},
	VOLUME = {11},
	YEAR = {2024},
	PAGES = {624--637},
	MRCLASS = {51F30 (05C05 53C23 57M05)},
	MRNUMBER = {4840249},
	DOI = {10.1090/bproc/243},
	URL = {https://doi.org/10.1090/bproc/243},
}

@article{saul2024some,
	title={Some novel constructions of optimal Gromov--Hausdorff--optimal correspondences between spheres},
	author={Rodr{\'i}guez Mart{\'i}n, Sa{\'u}l},
	journal={arXiv preprint arXiv:2409.02248},
	url={https://arxiv.org/abs/2409.02248},
	year={2024}
}

@article{memoli2025persistence,
  title={Persistence and Topological Complexity},
  author={Mémoli, Facundo and Zhou, Ling},
  journal={arXiv preprint arXiv:2506.17888},
  url={https://arxiv.org/abs/2506.17888},
  year={2025}
}

@article{memoli2012some,
	title     = {Some Properties of Gromov–Hausdorff Distances},
	author    = {Mémoli, Facundo},
	journal   = {Discrete \& Computational Geometry},
	volume    = {48},
	number    = {2},
	pages     = {416--440},
	year      = {2012},
	doi       = {10.1007/s00454-012-9406-8},
	url       = {https://doi.org/10.1007/s00454-012-9406-8}
}

@article{talipov2022gromov,
	AUTHOR = {Talipov, Talant K.},
	TITLE = {The {G}romov-{H}ausdorff distance between vertex sets of
	regular polygons inscribed in a single circle},
	NOTE = {Translation of Vestnik Moskov. Univ. Ser. I Mat. Mekh.
	{{\bf{2}}023}, no. 3, 23--27},
	JOURNAL = {Moscow Univ. Math. Bull.},
	FJOURNAL = {Moscow University Mathematics Bulletin},
	VOLUME = {78},
	YEAR = {2023},
	NUMBER = {3},
	PAGES = {130--135},
	ISSN = {0027-1322},
	MRCLASS = {54E35},
	MRNUMBER = {4641158},
	URL = {https://arxiv.org/abs/2210.09971}
}

@article{gakhar2019k,
	title={Künneth formulae in persistent homology},
	author={Gakhar, Hitesh and Perea, Jose A},
	eprint={1910.05656},
	journal={arXiv preprint arXiv:1910.05656},
	url={https://arxiv.org/abs/1910.05656},
	year={2019}
}

@article{zhou2023persistent,
	title={Persistent Sullivan minimal models of metric spaces},
	author={Zhou, Ling},
	eprint = {2310.06263},
	journal={arXiv preprint arXiv:2310.06263},
	url ={https://arxiv.org/abs/2310.06263},
	year={2023}
}

@article{memoli2024persistenthomotopy,
	title={Persistent homotopy groups of metric spaces},
	author={M{\'e}moli, Facundo and Zhou, Ling},
	journal={Journal of Topology and Analysis},
	pages={1--62},
	year={2024},
	publisher={World Scientific},
	doi       = {10.1142/S1793525324500018},
	url       = {https://doi.org/10.1142/S1793525324500018}
}

@article{memoli2025ephemeral,
	author    = {Facundo Mémoli and Ling Zhou},
	title     = {Ephemeral persistence features and the stability of filtered chain complexes},
	journal   = {Journal of Computational Geometry},
	volume    = {15},
	number    = {2},
	year      = {2024},
	pages     = {},
	doi       = {10.20382/jocg.v15.i2.a8},
	url       = {https://jocg.org/index.php/jocg/article/view/5193},
	note      = {Special Issue of Selected Papers from SoCG 2023}
}

@article{postol2023persistence,
	title={Algebraic Topology for Data Scientists},
	author={Michael S. Postol},
	year={2023},
	eprint={2308.10825},
	journal={arXiv preprint arXiv:2308.10825},
	primaryClass={math.AT},
	url={https://arxiv.org/abs/2308.10825},
}

@article{hess2024minimalmodels,
	title={Cell decompositions of persistent minimal models},
	author={Kathryn Hess and Samuel Lavenir and Kelly Maggs},
	year={2024},
	eprint={2312.08326},
	journal={arXiv preprint arXiv:2312.08326},
	primaryClass={math.AT},
	url={https://arxiv.org/abs/2312.08326},
}

@article{ginot2019distances,
	author       = {Gr\'egory Ginot and Johan Leray},
	title        = {Multiplicative persistent distances},
	year         = {2019},
	journal={arXiv preprint arXiv:1905.12307},
	url={https://arxiv.org/abs/1905.12307}
}

@article {belchi2022persistence,
	AUTHOR = {Belch\'i, Francisco and Stefanou, Anastasios},
	TITLE = {{$A_\infty$} persistent homology estimates detailed topology from pointcloud datasets},
	JOURNAL = {Discrete Comput. Geom.},
	FJOURNAL = {Discrete & Computational Geometry. An International Journal	of Mathematics and Computer Science},
	VOLUME = {68},
	YEAR = {2022},
	NUMBER = {1},
	PAGES = {274--297},
	ISSN = {0179-5376},
	MRCLASS = {55N31},
	MRNUMBER = {4430289},
	MRREVIEWER = {Marc Ethier},
	DOI = {10.1007/s00454-021-00319-y},
	URL = {https://doi.org/10.1007/s00454-021-00319-y},
}

@book {polterovich2020persistence,
	AUTHOR = {Polterovich, Leonid and Rosen, Daniel and Samvelyan, Karina
	and Zhang, Jun},
	TITLE = {Topological persistence in geometry and analysis},
	SERIES = {University Lecture Series},
	VOLUME = {74},
	PUBLISHER = {American Mathematical Society, Providence, RI},
	YEAR = {2020},
	PAGES = {xi+128},
	ISBN = {978-1-4704-5495-1},
	MRCLASS = {55N31 (53Dxx 58Cxx)},
	MRNUMBER = {4249570},
	MRREVIEWER = {Walter D. Freyn},
	URL = {https://arxiv.org/abs/1904.04044}
}

@article{carlsson2013viral,
	author    = {Chan, Joseph Minhow and Carlsson, Gunnar and Rabadan, Raul},
	title     = {Topology of viral evolution},
	journal   = {Proceedings of the National Academy of Sciences},
	volume    = {110},
	number    = {46},
	pages     = {18566--18571},
	year      = {2013},
	doi       = {10.1073/pnas.1313480110},
	url       = {https://doi.org/10.1073/pnas.1313480110},
}

@incollection {hausmann1995vietoris,
	AUTHOR = {Hausmann, Jean-Claude},
	TITLE = {On the {V}ietoris-{R}ips complexes and a cohomology theory for
	metric spaces},
	BOOKTITLE = {Prospects in topology ({P}rinceton, {NJ}, 1994)},
	SERIES = {Ann. of Math. Stud.},
	VOLUME = {138},
	PAGES = {175--188},
	PUBLISHER = {Princeton Univ. Press, Princeton, NJ},
	YEAR = {1995},
	MRCLASS = {57M99 (54E99 55N35)},
	MRNUMBER = {1368659},
	URL = {https://publish.illinois.edu/ymb/files/2020/03/Hausmann-1995-On-the-Vietoris-Rips-complexes-and-a-cohomology-th.pdf}
}

@article {gromov1983filling,
	AUTHOR = {Gromov, Mikhael},
	TITLE = {Filling {R}iemannian manifolds},
	JOURNAL = {J. Differential Geom.},
	FJOURNAL = {Journal of Differential Geometry},
	VOLUME = {18},
	YEAR = {1983},
	NUMBER = {1},
	PAGES = {1--147},
	ISSN = {0022-040X},
	MRCLASS = {53C20 (53C21 57R99)},
	MRNUMBER = {697984},
	MRREVIEWER = {Yu. Burago},
	URL = {http://projecteuclid.org/euclid.jdg/1214509283},
}

@phdthesis{aubrey2011thesis,
	AUTHOR = {HB, Aubrey},
	TITLE = {Persistent {C}ohomology {O}perations},
	NOTE = {Thesis (Ph.D.)--Duke University},
	PUBLISHER = {ProQuest LLC, Ann Arbor, MI},
	YEAR = {2011},
	PAGES = {119},
	ISBN = {978-1124-60902-7},
	MRCLASS = {Thesis},
	MRNUMBER = {2873406},
	URL = {http://gateway.proquest.com/openurl?url_ver=Z39.88-2004&rft_val_fmt=info:ofi/fmt:kev:mtx:dissertation&res_dat=xri:pqdiss&rft_dat=xri:pqdiss:3453041},
}

@article {bubenik2015metrics,
	AUTHOR = {Bubenik, Peter and de Silva, Vin and Scott, Jonathan},
	TITLE = {Metrics for generalized persistence modules},
	JOURNAL = {Found. Comput. Math.},
	FJOURNAL = {Foundations of Computational Mathematics. The Journal of the
	Society for the Foundations of Computational Mathematics},
	VOLUME = {15},
	YEAR = {2015},
	NUMBER = {6},
	PAGES = {1501--1531},
	ISSN = {1615-3375},
	MRCLASS = {55N35 (55U10)},
	MRNUMBER = {3413628},
	MRREVIEWER = {Mikael Vejdemo Johansson},
	DOI = {10.1007/s10208-014-9229-5},
	URL = {https://doi.org/10.1007/s10208-014-9229-5},
}

@article {blumberg2023interleaving,
	AUTHOR = {Blumberg, Andrew J. and Lesnick, Michael},
	TITLE = {Universality of the homotopy interleaving distance},
	JOURNAL = {Trans. Amer. Math. Soc.},
	FJOURNAL = {Transactions of the American Mathematical Society},
	VOLUME = {376},
	YEAR = {2023},
	NUMBER = {12},
	PAGES = {8269--8307},
	ISSN = {0002-9947},
	MRCLASS = {55P99 (55U99)},
	MRNUMBER = {4669297},
	DOI = {10.1090/tran/8738},
	URL = {https://doi.org/10.1090/tran/8738},
}

@article {chazal2014geometric,
	AUTHOR = {Chazal, Fr\'{e}d\'{e}ric and de Silva, Vin and Oudot, Steve},
	TITLE = {Persistence stability for geometric complexes},
	JOURNAL = {Geom. Dedicata},
	FJOURNAL = {Geometriae Dedicata},
	VOLUME = {173},
	YEAR = {2014},
	PAGES = {193--214},
	ISSN = {0046-5755},
	MRCLASS = {55N35 (55U10)},
	MRNUMBER = {3275299},
	MRREVIEWER = {Peter Bubenik},
	DOI = {10.1007/s10711-013-9937-z},
	URL = {https://doi.org/10.1007/s10711-013-9937-z},
}

@book {mosheroperations1968,
	AUTHOR = {Mosher, Robert E. and Tangora, Martin C.},
	TITLE = {Cohomology operations and applications in homotopy theory},
	PUBLISHER = {Harper \& Row, Publishers, New York-London},
	YEAR = {1968},
	PAGES = {x+214},
	MRCLASS = {55.40},
	MRNUMBER = {226634},
	MRREVIEWER = {Larry Smith},
	URL = {https://archive.org/details/cohomologyoperat0000mosh}
}

@article {Crawley-Boevey.2015,
	AUTHOR = {Crawley-Boevey, William},
	TITLE = {Decomposition of pointwise finite-dimensional persistence modules},
	JOURNAL = {J. Algebra Appl.},
	FJOURNAL = {Journal of Algebra and its Applications},
	VOLUME = {14},
	YEAR = {2015},
	NUMBER = {5},
	PAGES = {1550066, 8},
	ISSN = {0219-4988},
	MRCLASS = {16G20},
	MRNUMBER = {3323327},
	MRREVIEWER = {M\'{a}ty\'{a}s Domokos},
	DOI = {10.1142/S0219498815500668},
	URL = {https://doi.org/10.1142/S0219498815500668},
}

@article {Chazal.2016b,
	AUTHOR = {Chazal, Fr\'{e}d\'{e}ric and Crawley-Boevey, William and de Silva,
	Vin},
	TITLE = {The observable structure of persistence modules},
	JOURNAL = {Homology Homotopy Appl.},
	FJOURNAL = {Homology, Homotopy and Applications},
	VOLUME = {18},
	YEAR = {2016},
	NUMBER = {2},
	PAGES = {247--265},
	ISSN = {1532-0073},
	MRCLASS = {55U05 (55N99)},
	MRNUMBER = {3575998},
	MRREVIEWER = {Ne\v{z}a Mramor Kosta},
	DOI = {10.4310/HHA.2016.v18.n2.a14},
	URL = {https://doi.org/10.4310/HHA.2016.v18.n2.a14},
}

@book {Chazal.2016a,
	AUTHOR = {Chazal, Fr\'{e}d\'{e}ric and de Silva, Vin and Glisse, Marc and Oudot,
	Steve},
	TITLE = {The structure and stability of persistence modules},
	SERIES = {Springer Briefs in Mathematics},
	PUBLISHER = {Springer},
	YEAR = {2016},
	PAGES = {x+120},
	ISBN = {978-3-319-42543-6},
	MRCLASS = {55N10 (16G20 55U10)},
	MRNUMBER = {3524869},
	MRREVIEWER = {Henry Hugh Adams},
	DOI = {10.1007/978-3-319-42545-0},
	URL = {https://doi.org/10.1007/978-3-319-42545-0},
}

@article {azumaya1950theorem,
	AUTHOR = {Azumaya, Gor\^{o}},
	TITLE = {Corrections and supplementaries to my paper concerning
	{K}rull-{R}emak-{S}chmidt's theorem},
	JOURNAL = {Nagoya Math. J.},
	FJOURNAL = {Nagoya Mathematical Journal},
	VOLUME = {1},
	YEAR = {1950},
	PAGES = {117--124},
	ISSN = {0027-7630},
	MRCLASS = {09.1X},
	MRNUMBER = {37832},
	MRREVIEWER = {I. Kaplansky},
	URL = {http://projecteuclid.org/euclid.nmj/1118764711},
}

@article{adamaszek2020homotopy,
	AUTHOR = {Adamaszek, Micha{\l} and Adams, Henry and Gasparovic, Ellen and Gommel, Maria and Purvine, Emilie and Sazdanovic, Radmila and Wang, Bei and Wang, Yusu and Ziegelmeier, Lori},
	TITLE = {On homotopy types of {V}ietoris--{R}ips complexes of metric gluings},
	JOURNAL = {J. Appl. Comput. Topol.},
	FJOURNAL = {Journal of Applied and Computational Topology},
	VOLUME = {4},
	YEAR = {2020},
	NUMBER = {3},
	PAGES = {425--454},
	ISSN = {2367-1726},
	MRCLASS = {55N31 (55U10)},
	MRNUMBER = {4130978},
	MRREVIEWER = {Jerzy Jezierski},
	DOI = {10.1007/s41468-020-00054-y},
	URL = {https://doi.org/10.1007/s41468-020-00054-y},
}

@book{burago2001course,
	AUTHOR = {Burago, Dmitri and Burago, Yuri and Ivanov, Sergei},
	TITLE = {A course in metric geometry},
	SERIES = {Graduate Studies in Mathematics},
	VOLUME = {33},
	PUBLISHER = {American Mathematical Society, Providence, RI},
	YEAR = {2001},
	PAGES = {xiv+415},
	ISBN = {0-8218-2129-6},
	MRCLASS = {53C23},
	MRNUMBER = {1835418},
	MRREVIEWER = {Mario Bonk},
	DOI = {10.1090/gsm/033},
	URL = {https://doi.org/10.1090/gsm/033},
}

@article{adams2022metric,
	AUTHOR = {Adams, Henry and Heim, Mark and Peterson, Chris},
	TITLE = {Metric thickenings and group actions},
	JOURNAL = {J. Topol. Anal.},
	FJOURNAL = {Journal of Topology and Analysis},
	VOLUME = {14},
	YEAR = {2022},
	NUMBER = {3},
	PAGES = {587--613},
	ISSN = {1793-5253},
	MRCLASS = {55N31 (05E45 20F65 54E35 55P10)},
	MRNUMBER = {4493474},
	MRREVIEWER = {Walter D. Freyn},
	DOI = {10.1142/S1793525320500569},
	URL = {https://doi.org/10.1142/S1793525320500569},
}

@article{lim2024vietoris,
	AUTHOR = {Lim, Sunhyuk and M\'{e}moli, Facundo and Okutan, Osman Berat},
	TITLE = {Vietoris--{R}ips persistent homology, injective metric spaces,
	and the filling radius},
	JOURNAL = {Algebr. Geom. Topol.},
	FJOURNAL = {Algebraic \& Geometric Topology},
	VOLUME = {24},
	YEAR = {2024},
	NUMBER = {2},
	PAGES = {1019--1100},
	ISSN = {1472-2747},
	MRCLASS = {55N31 (53C23)},
	MRNUMBER = {4735051},
	DOI = {10.2140/agt.2024.24.1019},
	URL = {https://doi.org/10.2140/agt.2024.24.1019},
}

@article{katz1983filling,
	AUTHOR = {Katz, Mikhail},
	TITLE = {The filling radius of two--point homogeneous spaces},
	JOURNAL = {J. Differential Geom.},
	FJOURNAL = {Journal of Differential Geometry},
	VOLUME = {18},
	YEAR = {1983},
	NUMBER = {3},
	PAGES = {505--511},
	ISSN = {0022-040X},
	MRCLASS = {53C20 (53C35)},
	MRNUMBER = {723814},
	MRREVIEWER = {Rolf Sulanke},
	URL = {http://projecteuclid.org/euclid.jdg/1214437785},
}

@book{hatcher2000,
	AUTHOR = {Hatcher, Allen},
	TITLE = {Algebraic topology},
	PUBLISHER = {Cambridge University Press, Cambridge},
	YEAR = {2002},
	PAGES = {xii+544},
	MRCLASS = {55-01 (55-00)},
	MRNUMBER = {1867354},
	MRREVIEWER = {Donald W. Kahn},
	URL = "https://cds.cern.ch/record/478079"
}

@article{memoli2024persistent,
	AUTHOR = {M\'{e}moli, Facundo and Stefanou, Anastasios and Zhou, Ling},
	TITLE = {Persistent cup product structures and related invariants},
	JOURNAL = {J. Appl. Comput. Topol.},
	FJOURNAL = {Journal of Applied and Computational Topology},
	VOLUME = {8},
	YEAR = {2024},
	NUMBER = {1},
	PAGES = {93--148},
	ISSN = {2367-1726},
	MRCLASS = {55U99 (55M30 55N20)},
	MRNUMBER = {4707458},
	DOI = {10.1007/s41468-023-00138-5},
	URL = {https://doi.org/10.1007/s41468-023-00138-5},
}

@article{adamaszek2018metric,
	AUTHOR = {Adamaszek, Micha{\l} and Adams, Henry and Frick, Florian},
	TITLE = {Metric reconstruction via optimal transport},
	JOURNAL = {SIAM J. Appl. Algebra Geom.},
	FJOURNAL = {SIAM Journal on Applied Algebra and Geometry},
	VOLUME = {2},
	YEAR = {2018},
	NUMBER = {4},
	PAGES = {597--619},
	MRCLASS = {54E35 (53C23 55P10 55U10)},
	MRNUMBER = {3871057},
	MRREVIEWER = {Thomas H\"{u}ttemann},
	DOI = {10.1137/17M1148025},
	URL = {https://doi.org/10.1137/17M1148025},
}

@article{gillespie2024vietoris,
	AUTHOR = {Gillespie, Patrick},
	TITLE = {Vietoris thickenings and complexes are weakly homotopy equivalent},
	JOURNAL = {J. Appl. Comput. Topol.},
	FJOURNAL = {Journal of Applied and Computational Topology},
	VOLUME = {8},
	YEAR = {2024},
	NUMBER = {1},
	PAGES = {35--53},
	ISSN = {2367-1726},
	MRCLASS = {55N31 (55P10)},
	MRNUMBER = {4707456},
	DOI = {10.1007/s41468-023-00135-8},
	URL = {https://doi.org/10.1007/s41468-023-00135-8},
}

@article{lovasz1983self,
	AUTHOR = {Lov{\'a}sz, L{\'a}sl{\'o}},
	TITLE = {Self--dual polytopes and the chromatic number of distance graphs on the sphere},
	JOURNAL = {Acta Sci. Math. (Szeged)},
	FJOURNAL = {Acta Universitatis Szegediensis. Acta Scientiarum
	Mathematicarum},
	VOLUME = {45},
	YEAR = {1983},
	NUMBER = {1-4},
	PAGES = {317--323},
	ISSN = {0001-6969},
	MRCLASS = {05C15 (52A25)},
	MRNUMBER = {708798},
	MRREVIEWER = {J. Spencer},
	URL = {http://pub.acta.hu/acta/showCustomerArticle.action?id=12377&dataObjectType=article&returnAction=showCustomerVolume&sessionDataSetId=7bedb62ae745f2a0&style=}
}

@incollection{contessoto_et_al:LIPIcs.SoCG.2022.31,
	AUTHOR = {Contessoto, Marco and M\'{e}moli, Facundo and Stefanou, Anastasios and Zhou, Ling},
	TITLE = {Persistent cup-length},
	BOOKTITLE = {38th {I}nternational {S}ymposium on {C}omputational {G}eometry},
	%SERIES = {LIPIcs. Leibniz Int. Proc. Inform.},
	VOLUME = {224},
	PAGES = {Art. No. 31, 17},
	PUBLISHER = {Schloss Dagstuhl. Leibniz-Zent. Inform., Wadern},
	YEAR = {2022},
	MRCLASS = {55N31 (68U05)},
	MRNUMBER = {4470910},
	DOI = {10.4230/lipics.socg.2022.31},
	URL = {https://doi.org/10.4230/lipi	cs.socg.2022.31},
}

@phdthesis{zhou2023beyond,
	title={Beyond Persistent Homology: More Discriminative Persistent Invariants},
	author={Zhou, Ling},
	year={2023},
	school={The Ohio State University},
	url={http://rave.ohiolink.edu/etdc/view?acc_num=osu1689837764936381}
}

@article{huang2005cup,
	title={Cup products in computational topology},
	author={Huang, Jonathan},
	year={2005},
	publisher={Citeseer},
	url={https://citeseerx.ist.psu.edu/document?repid=rep1&type=pdf&doi=04b5ed449bff26197856c8ecb7dd4fbbf93ea846}
}

@phdthesis{yarmola2010persistence,
	title={Persistence and computation of the cup product},
	author={Yarmola, Andrew},
	advisor = {Mikael Vejdemo-Johansson},
	year={2010},
	school={Stanford University},
	type={Undergraduate honors thesis},
	URL = {https://web.math.princeton.edu/~yarmola/assets/pdf/cup_prod.pdf}
}

@article{herscovich2018higher,
	AUTHOR = {Herscovich, Estanislao},
	TITLE = {A higher homotopic extension of persistent (co)homology},
	JOURNAL = {J. Homotopy Relat. Struct.},
	FJOURNAL = {Journal of Homotopy and Related Structures},
	VOLUME = {13},
	YEAR = {2018},
	NUMBER = {3},
	PAGES = {599--633},
	ISSN = {2193-8407},
	MRCLASS = {55N35 (16E45 16W70 18G55 55U10)},
	MRNUMBER = {3856303},
	MRREVIEWER = {Timothy Porter},
	DOI = {10.1007/s40062-017-0195-x},
	URL = {https://doi.org/10.1007/s40062-017-0195-x},
}

@phdthesis{contreras2022persistent,
	AUTHOR = {Contreras, Luis G. Polanco},
	TITLE = {Applications of {P}ersistent {C}ohomology to {D}imensionality {R}eduction and {C}lassification {P}roblems},
	NOTE = {Thesis (Ph.D.)--Michigan State University},
	PUBLISHER = {ProQuest LLC, Ann Arbor, MI},
	YEAR = {2022},
	PAGES = {88},
	MRCLASS = {Thesis},
	MRNUMBER = {4435582},
	URL = {http://gateway.proquest.com/openurl?url_ver=Z39.88-2004&rft_val_fmt=info:ofi/fmt:kev:mtx:dissertation&res_dat=xri:pqm&rft_dat=xri:pqdiss:29206251},
}

@unpublished{barham2025group,
	author = {Barham, Liam},
	title = {Group Actions on Metric Spaces},
	note = {Talk presented at the 58th Spring Topology and Dynamics Conference},
	date = {2025-03-08},
	url = {https://scholarlattice.org/events/a61e84e5-8074-4f9c-8f7d-7f5479164a14},
	urldate = {2026-09-01},
}

@inproceedings{chazal2009gromov,
	title={Gromov--Hausdorff stable signatures for shapes using persistence},
	author={Chazal, Fr{\'e}d{\'e}ric and Cohen-Steiner, David and Guibas, Leonidas J and M{\'e}moli, Facundo and Oudot, Steve Y},
	booktitle={Computer Graphics Forum},
	volume={28},
	number={5},
	pages={1393--1403},
	year={2009},
	doi = {10.1111/j.1467-8659.2009.01516.x},
	url = {https://doi.org/10.1111/j.1467-8659.2009.01516.x}
}

@article {adamaszek2017VietorisProduct,
	AUTHOR = {Adamaszek, Micha{\l} and Adams, Henry},
	TITLE = {The {V}ietoris-{R}ips complexes of a circle},
	JOURNAL = {Pacific J. Math.},
	FJOURNAL = {Pacific Journal of Mathematics},
	VOLUME = {290},
	YEAR = {2017},
	NUMBER = {1},
	PAGES = {1--40},
	ISSN = {0030-8730},
	MRCLASS = {55P15 (05C20 05E45 55U10)},
	MRNUMBER = {3673078},
	MRREVIEWER = {Thomas H\"{u}ttemann},
	DOI = {10.2140/pjm.2017.290.1},
	URL = {https://doi.org/10.2140/pjm.2017.290.1},
}
\end{document}